\documentclass[11pt,reqno]{amsart}
\usepackage{graphicx} % Required for inserting images
\usepackage{comment}
\usepackage{amsthm}
\usepackage{amsmath}
\usepackage{amsfonts}
\usepackage{amssymb}
\usepackage{color}
\usepackage{graphicx}
\usepackage{url}
\usepackage{hyperref}
\usepackage{comment}

\usepackage[normalem]{ulem}

\newtheorem{theorem}{Theorem}[section]

\newtheorem{proposition}[theorem]{Proposition}

\theoremstyle{definition}
\newtheorem{definition}[theorem]{Definition}
\theoremstyle{remark}
\newtheorem{remark}[theorem]{Remark}
\newtheorem{example}[theorem]{Example}

\numberwithin{equation}{section}

\allowdisplaybreaks

\title[The $p$-Laplacian: flow in porous media and CFD simulations]{The $p$-Laplacian: phenomenological modelling of the flow in porous media and CFD simulations}
\author[P.~Girg]{Petr Girg}
\author[L.~Kotrla]{Lukáš Kotrla}
\author[A.~Švandová]{Anežka Švandová}

\address{Petr Girg \newline
Department of Mathematics and NTIS,
Faculty of Applied Scences, University of West Bohemia,
Univerzitn\'{\i} 8, CZ-301\,00~Plze\v{n}, Czech Republic}
\email{pgirg@kma.zcu.cz}

\address{Luk\'a\v{s} Kotrla \newline
Department of Mathematics and NTIS,
Faculty of Applied Scences, University of West Bohemia,
Univerzitn\'{\i} 8, CZ-301\,00~Plze\v{n}, Czech Republic}
\email{kotrla@ntis.zcu.cz}

\address{Ane\v{z}ka \v{S}vandov\'{a} \newline
Department of Mathematics,
Faculty of Applied Scences, University of West Bohemia,
Univerzitn\'{\i} 8, CZ-301\,00~Plze\v{n}, Czech Republic}
\email{svandova@kma.zcu.cz}

\begin{document}

\keywords{$p$-Laplacian; Porous medium; filtration; Darcy's law;
\hfill\break\indent pressure-to-velocity power law}

\begin{abstract}
The aim of this paper is to discuss several aspects of connections between the
$p$-Laplacian and mathematical models in hydrology. 
At first we present models of groundwater flow in phreatic aquifers and models of irrigation and drainage that lead to quasilinear parabolic equations involving the $p$-Laplacian. 
Next, we survey conditions of validity of Strong Maximum Principle and Strong Comparison Principle for this type of problems. Finally, we employ computer fluid dynamics simulations to realistic scenario of fracture networks
to estimate values of the parameters
of constitutive laws governing groundwater flow in the context of fractured
hard-rock aquifers.   
\end{abstract}

\maketitle

\section{Introduction}

The aim of this paper is to discuss several aspects of connections between 
the $p$-Laplacian and mathematical models of groundwater flow with applications to irrigation, drainage, and fresh water supply in small rural areas.
The operator $p$-Laplacian, a non-linear generalization of the Laplace operator,
$u\mapsto \mathop{\mathrm{div}}\left(|\mathop{\nabla}u|^{p-2}\mathop{\nabla }u\right)$, $p>1, p\not=2$, gained substantial attention of mathematicians working in nonlinear functional analysis since the 1970s.
This can be attributed mostly to the fact that the $p$-Laplacian exhibits many features of more general nonlinear operators, while being simple and elegant, and hence allowing the ideas of the proofs to be kept clear and accessible. However, origins of the $p$-Laplacian can be traced back to the 1870s, see, e.g.,  \cite{Smreker1878, Smreker1879, Zhukovskii1889}, when the radially symmetric version of this operator (written in a different way than it is customary today) was used in the theoretical study of groundwater flow towards a well in coarse grained porous media (such as gravels). To the best of our knowledge, the $p$-Laplacian in its full PDE form was first introduced in the paper \cite{Leibenson1945}
in the context of mathematical models of the flow of natural gas in the porous rock
as early as in 1945. In particular, the following equation
(written in modern notation)
$$
u^{-\kappa}(\boldsymbol{x},t)\frac{\partial u}{\partial t}(\boldsymbol{x},t)=
c\mathop{\Delta_{p}}u(\boldsymbol{x},t)\,,
$$ where $c>0$, $\kappa=1/2$, $p=3/2$, $\boldsymbol{x}\in\mathbb{R}^3$, and $t\in (0, +\infty)$,
was proposed as a model of isothermic process of turbulent filtration of natural gas in a porous medium in \cite{Leibenson1945}. This equation was soon after suggested in a more general setting as a model of polytropic process of turbulent filtration of natural gas with $\kappa=\gamma/(\gamma+1)$, where $\gamma$ is the polytropic index of the process and $p\in [3/2, 2]$, see \cite[p.~504]{Leibenson1953}.
In the early papers on the subject, authors (without the benefit of
techniques of modern functional analysis) limited themselves to study particular cases of initial-boundary value problems with the $p$-Laplacian which were important for engineers of that days, with breakthrough results in \cite{Barenblatt1952b, Barenblatt1952}. More information about pioneering works related to applications of the $p$-Laplacian to flow in porous media can be found in \cite{BenediktGirgKotrlaTakac2018}. 

In this paper, we cover two topics connecting the $p$-Laplacian and mathematical models in hydrology. The first topic is related to mathematical models in irrigation and drainage and concerns validity of Strong Maximum Principle (SMP) and Strong Comparison Principle (SCP) for initial-value problem
\begin{equation}
\label{in:prob}
\left\{\begin{aligned}
\frac{\partial}{\partial t} b(u(\boldsymbol{x}, t))-c\,\Delta_p u(\boldsymbol{x}, t) & =f(\boldsymbol{x}, t) \geq 0 & & \text { for }(\boldsymbol{x}, t) \in \Omega \times(0, T) ; \\
u(\boldsymbol{x}, 0) & =u_0(\boldsymbol{x}) \geq 0 & & \text { for } \boldsymbol{x} \in \Omega ; \\
u(\boldsymbol{x}, t) & =0 & & \text { for }(\boldsymbol{x}, t) \in \partial \Omega \times(0, T)\,, \\
u(\boldsymbol{x}, t) & \geq 0 & & \text { for }(\boldsymbol{x}, t) \in \overline{\Omega} \times[0, T]\,,
\end{aligned}\right.
\end{equation}
where $\Omega\subset\mathbb{R}^N$, $N\in\mathbb{N}$, is a bounded domain, $c\equiv\mathrm{const.}>0$, 
$b: \mathbb{R}_{+} \rightarrow \mathbb{R}_{+}:=[0,+\infty)$ is a continuous function, $b(0)=0$, and $b \in C^1(0,+\infty)$ with $b^{\prime}>0$ in $(0,+\infty)$,
$f\in C(\Omega \times(0, T))$ and $u_0\in C(\Omega)$ 
are nonnegative functions. 
For detailed discussion about mathematical models used in irrigation and drainage based on~\eqref{in:prob}, 
see Section~2.
Explicit solutions of problem~\eqref{in:prob} can be obtained only in rare special cases. Thus we rely on qualitative analysis and numerical methods in dealing with \eqref{in:prob}. 
Maximum and comparison principles play important role in this process.  
In this paper, we discuss conditions under which these principles
for \eqref{in:prob} hold and provide some realistic examples when they do not hold. This is not interesting only from theoretical point of view, but it has implications for the choices of appropriate numerical methods to find numerical solutions of~\eqref{in:prob}. 

The second topic delves into the utilization of
computer fluid dynamics (CFD) numerical simulations to estimate values of the parameters of constitutive laws governing groundwater flow in the context of fractured hard rock. 
This research builds
upon existing research, see, e.g.,~\cite{Brush2003, sarkar2004fluid, StarkVolker,  yan2018non, Zimmerman2004163},
and is motivated by the fact that increasing demand for fresh water has driven interest in hard-rock aquifers \cite{Gustafson1994, Shapiro2002}, despite their limited well yields due to water flow occurring only in cracks and fractures. With 
hard-rocks covering over 20\%
of the Earth's landmass \cite{AYRAUD20082686, Lachassagne2021, Gustafson1994},
these aquifers hold significant freshwater reserves, particularly in semi-arid regions like sub-Saharan Africa \cite{ MacDonald2008,Wright1992}, Australia \cite{Australia2014} and India \cite{Perrin2011, Singhal2008}.
They provide important water sources for rural populations in these areas, particularly in Australia  and India, which store 40\% and 50\% of their groundwater in such aquifers, respectively.
Understanding flow dynamics in these aquifers is essential for improving rural living standards through better access to fresh water.

The paper is organized as follows. In Section~\ref{sec:mathmod}, we present several mathematical models of groundwater flow in phreatic aquifers and related models used in irrigation and drainage. In Section~\ref{sec:maxprincip},
we survey some of our recent results concerning SMP and SCP.
In Section~\ref{sec:CFD}, we present some of our recent results concerning estimations of parameters
in constitutive laws
governing groundwater flow through fractured hard rocks.
Section~\ref{sec:drabek} is devoted to contribution of P.~Dr\'{a}bek to topics related to the $p$-Laplacian.
This paper adheres to the SI system for all physical quantities (m for length, s for time, kg for mass, etc.).

\section{Mathematical
models}
\label{sec:mathmod}

In mathematical models of groundwater flow, averaged velocity, total head, and piezometric head are crucial concepts.
Actual velocity of the groundwater highly oscillates in the channels in the porous medium and thus it is difficult to measure and predict.
Averaged velocity
$\vec{v}_{\mathrm{av}}$ captures the idea of bulk motion of the groundwater within a sufficiently large control volume in the porous medium.
It is defined by means of specific discharge $\vec{q}$ by formula $\vec{v}_{\mathrm{av}}=\vec{q}/\phi_{\mathrm{eff}}$, where the specific discharge (vector) $\vec{q}$ takes the direction of the flow and its magnitude is defined as the volume of water flowing per unit time through a unit cross-sectional area normal to the direction of flow, and $\phi_{\mathrm{eff}}$ is effective porosity of the medium
see, e.g., 
\cite[p.~121]{Bear1972} for detailed explanation.
On this macroscopic level, the process of transformation and dissipation of energy can be described in the following way. The total mechanical  energy per volume $E_T$ in a control volume of water is the sum of gravitational potential energy $z \varrho g$, pressure energy $P$, and kinetic energy $\frac{1}{2} \varrho {v}_{\mathrm{av}}^{2}$ (all three per volume), where ${v}_{\mathrm{av}}$ stands for the magnitude of averaged velocity of the flow in the control volume, $\varrho$ is the water density, $P$ pressure, $z$ elevation of the control volume from the datum, $g$ gravitational acceleration, see, e.g., \cite{ZekaiSen1995}. For incompressible liquid such as water, one can equivalently consider another quantity 
$$
h_{T}=\frac{E_T}{\varrho g}=z+\frac{P}{\varrho g}+\frac{1}{2 g} v_{\mathrm{av}}^{2}
$$
which can be directly measured in practice, e.g., by using observation wells or by so called piezometers, see, e.g., \cite[p.~63]{Bear1972} or \cite[pp.~129--132]{ZekaiSen1995}. Groundwater is loosing its total energy (or equivalently total head $h_T$) while flowing due to viscous forces and friction with porous medium. Thus, its total energy decreases in the direction of the flow. In typical real-world situations, the term corresponding to kinetic energy is negligible and can be dropped, see, e.g., 
\cite[p.~5]{Harr2012groundwater} or
\cite[pp.~40--43]{ZekaiSen1995}. 
In this way, we obtain piezometric head
$$
h=z+\frac{P}{\varrho g}
$$
which is the state variable in the mathematical models of the water flow in the underground. On the other hand, the specific discharge is the flux quantity.
The constitutive law relating this two quantities, quantitatively describes the rate of dissipation of the energy along the flow path.

In this paper, we will limit ourselves to a mathematical model of unconfined
aquifer bounded below by a flat impermeable layer at $z=0$. 
This constitutes a free boundary problem in its full generality, since the upper boundary is the unknown surface of the groundwater. In groundwater modelling, this difficulty is simplified by assuming that the vertical flux in the aquifer is negligible, leading to
piezometric head $h$ being constant in the $z$-direction. This assumption is known as Dupuit-Forchheimer assumption, see, e.g., \cite{Bear1972, Bear2014, ArielSpringer2021} and it is based on observations performed on aquifers.
Then the height $\widehat{h}$ of the free surface of the groundwater (called water table for short)  above the point $(x,y,0)$ at time $t$ is 
$\widehat{h}(x,y,t)=h(x,y,0,t)$. 
The balance equation can be written as follows
\begin{equation}
\label{eq:balance}    
\phi_{\mathrm{eff}} \frac{\partial\widehat{h}}{\partial t}(x,y,t) 
    - 
    \mathop{ \mathrm{div} } \left( \widehat{h}\left(x,y,t\right)\ {\vec{q}}\left(x,y,t\right) \right)
    = \widehat{g}(x,y,t)\,,
\end{equation}
see, e.g., \cite[Eq.~(5.4.43)]{Bear2014}, 
\cite[Eq.~2.6]{ArielSpringer2021}.
Here, $\widehat{g}(x,y,t)$ represents external sinks and sources (evaporation, rainfall etc.) of volume of water per area and time unit. Recall that 
$\phi_{\mathrm{eff}}$ is effective porosity of the porous medium
and $\vec{q}(x,y,t)$ is specific discharge.

To eliminate flux variable from the balance equation \eqref{eq:balance},
suitable constitutive law is needed.
It is obtained empirically from experimental data for given porous medium and fluid. In practice, these experiments are performed on 
a~sample of porous medium subjected to one dimensional flow for several values of magnitudes of flux $q$, which are kept constant during each measurement. 
Linear Darcy's law is the most widely used in practice due to its simplicity and still reasonable accuracy. It relates groundwater flux to the piezometric head loss per length according to the following formula
\begin{equation}
\label{const:Darcy}
q = c\, \frac{\triangle h}{\triangle L}\,,
\end{equation}
where $c>0$ is a constant to be determined from measured data,  
$\triangle h$ is the difference of the piezometric heads measured at two distinct locations distance $\triangle L$ apart. 
This formula was established experimentally for filtration of water through sand by Henry Darcy \cite{Darcy1856} in 1856.
Later, it was found that it has limited range of its validity in coarse grained media (such as gravels), see, e.g., \cite{Forchheimer1901, Izbash1931, Kroeber1884, Missbach1936, Smreker1878, Smreker1879} as well as in media with very low permeability (such as clays, certain soils, and sandstones), see, e.g., \cite{CHENG2024, King1898, Zunker1920}. For thorough surveys and discussions of this and other constitutive laws and various criteria of their validity, see, e.g., \cite{AravinNumerov, 
Bear1972, Scheidegger1960,
ZekaiSen1995,
SoniIslamBasak1978}. It follows from discussions in these papers that the 
power-type law
\begin{equation}
\label{power:law}
q = c\, \left(\frac{\triangle h}{\triangle L}\right)^m
\end{equation}
with  constants $c, m>0$ to be determined from measured data,
is simple but flexible enough to fit with most experimental data obtained for various porous media. 
For $m=1$, the power-type law coincides with Darcy's law. 
The case $m>1$
correponds to natural media such as clays, certain soils and sandstones, while the case
$1/2<m<1$
corresponds to coarse grained materials such as gravels, see the literature listed above. 
Let us note that the law~\eqref{power:law} is not the only type of nonlinear laws used in practice. For thorough surveys, see, e.g., \cite{AravinNumerov,BenediktGirgKotrlaTakac2018, Bear1972, Scheidegger1960, ZekaiSen1995}.

For  the homogeneous and isotropic porous medium, the two- or three-dimensional constitutive law in differential form can be inferred from the 
one-dimensional one, by taking into account that
the flux takes the opposite direction of the gradient of the piezometric head and no 
flow occurs if the gradient of the piezometric head is zero. In this way, we obtain
\begin{equation}
\label{vec:power:law}
\vec{q} = 
\begin{cases} \vec{0} & \text { for } \nabla h=\vec{0}\,, \\ 
- c\, \left|\nabla h\right|^{m-1} \nabla h  & 
\text { for } 
\nabla h\neq\vec{0}\,,
\end{cases}
\end{equation}
where $\nabla h$ stands for the spatial gradient of the piezometric head, $c, m>0$ are constants as in \eqref{power:law}.

By substituting two-dimensional differential version of Darcy's law, i.e., \eqref{vec:power:law} with $m=1$ for $\vec{q}$ in \eqref{eq:balance}, we obtain classical porous medium equation
$$
\phi_{\mathrm{eff}}\frac{\partial \widehat{h}}{\partial t}-
c\operatorname{div}(\widehat{h} \nabla \widehat{h})=\widehat{g}(x, y, t)\,.
$$
By the same procedure, we obtain the following nonlinear equation
\begin{equation}
\label{eq:doubly:hath}
\phi_{\mathrm{eff}}\frac{\partial \widehat{h}}{\partial t}-
c\operatorname{div}(\widehat{h} |\nabla \widehat{h}|^{m-1}\nabla \widehat{h})=\widehat{g}(x, y, t)\,,
\end{equation}
by using \eqref{vec:power:law} with $m>0$, $m\not=1$.

Now we turn our attention to mathematical models from irrigation and drainage. We were motivated by models proposed in \cite{SinghRai2D}, but we take into account nonlinear effects and use the
power-type law~\eqref{power:law} instead of Darcy's law \eqref{const:Darcy}. 

At first, let us consider local aquifer under a field, represented by a bounded domain $\Omega\subset\mathbb{R}^2$, surrounded by open water body as in Figure~1. The horizontal permeable layer is bounded from below by the horizontal impermeable layer (bedrock). The drainage channels are fully penetrating, i.e., they reach the surface of the impermeable 
bedrock.  We take the bedrock surface as vertical datum, i.e., we assign to it coordinate $z=0$. Let $H>0$ be the depth of the open water bodies relative to the datum. Then the height $\widehat{h}$ of the water table above the datum
satisfies balance equation \eqref{eq:doubly:hath} with boundary conditions $\widehat{h}=H$ on $\partial\Omega$. 
Using substitution $p=m+1$ (to match notation with the $p$-Laplacian),
$u=\widehat{h}^{p /(p-1)}-H^{p /(p-1)}$, we obtain initial-boundary value problem \eqref{in:prob} with
\begin{equation}
\label{eq:bu:water}
b(u)=
\phi_{\mathrm{eff}}\left(\frac{p}{p-1}\right)^{p-1}
\left[
\left(u+H^{\frac{p}{p-1}}\right)^{\frac{p-1}{p}}-H
\right]\,,
\end{equation} 
$f(\boldsymbol{x},t)=\left(\frac{p}{p-1}\right)^{p-1} \widehat{g}(\boldsymbol{x},t)$,
and $u_0(\boldsymbol{x})=\widehat{h}_0^{p /(p-1)}(\boldsymbol{x})-H^{p /(p-1)}$, where $\boldsymbol{x}=(x,y)\in\Omega$.

The second model covered by our theoretical methods concerns local phreatic aquifer under a strip field between two parallel ditches. The aquifer is bounded from below by an impermeable bedrock and the ditches are fully penetrating to the bedrock, see Figure~\ref{fig:parallel}. In general, the domain would be infinite strip in this situation. However, if we assume translation invariance of the solution along the parallel ditches, we can limit ourselves to a bounded interval $\Omega=(-L/2,L/2)$, where $L>0$, see Figure~2. This nonlinear model has been proposed in~\cite{ArielSpringer2021} and it is motivated by~~\cite{BASAK1979,Marino1974,SinghRai1980,SinghRai1989}, where linear Darcy's law~\eqref{const:Darcy} is used instead. With $\widehat{h}$ being piezometric head and the bedrock vertical datum, the function $u=\widehat{h}^{p /(p-1)}-H^{p /(p-1)}$ satisfies
\eqref{in:prob} with $b, f, u_0$ as above but with  $\boldsymbol{x}=x\in\Omega=(-L/2, L/2)$. 

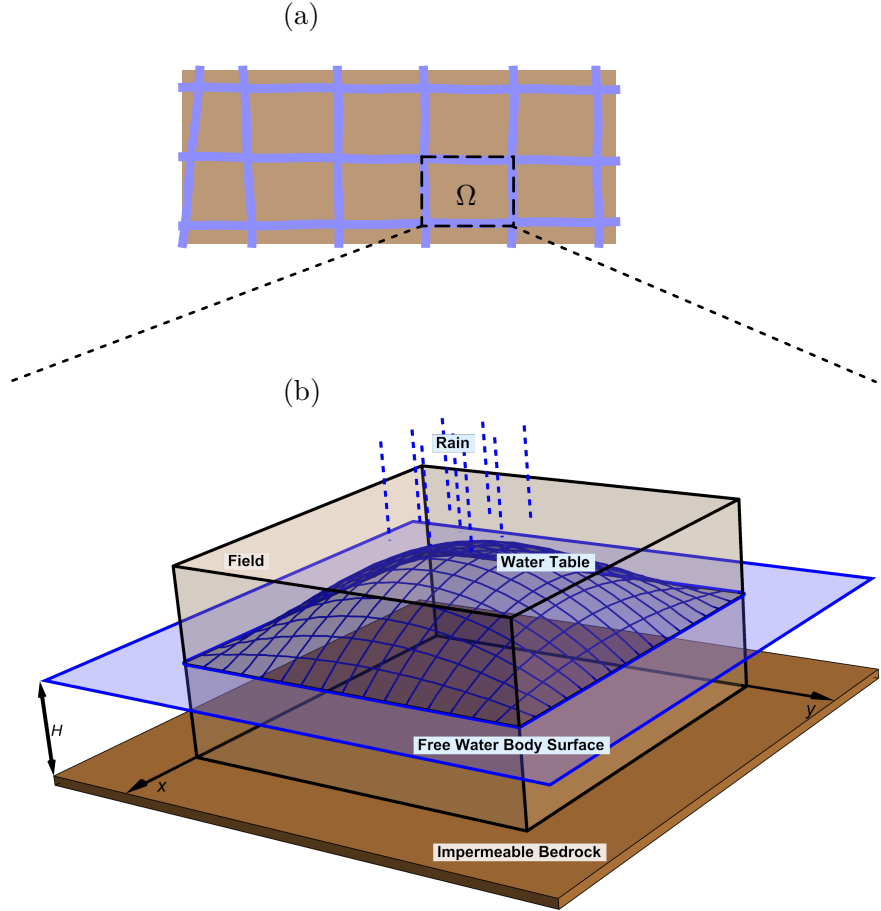
\begin{figure}[htp]
\label{fig:grid:fields}
\setlength{\unitlength}{1cm}
\begin{picture}(12.5,12.5)
\put(0,7){
\includegraphics[width=12cm]{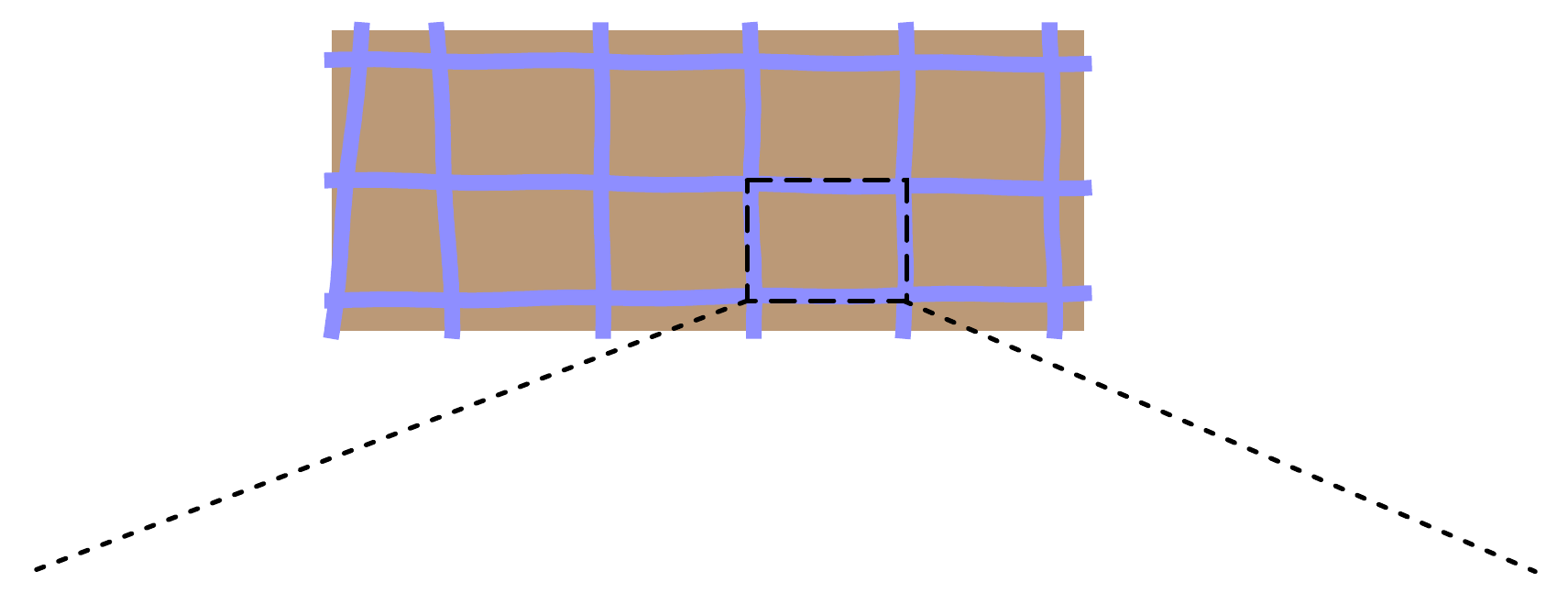}}
\put(4,12){(a)}
\put(6.3,9.6){$\Omega$}
\put(4,7){(b)}
\put(0,0){
\includegraphics[width=12cm]
{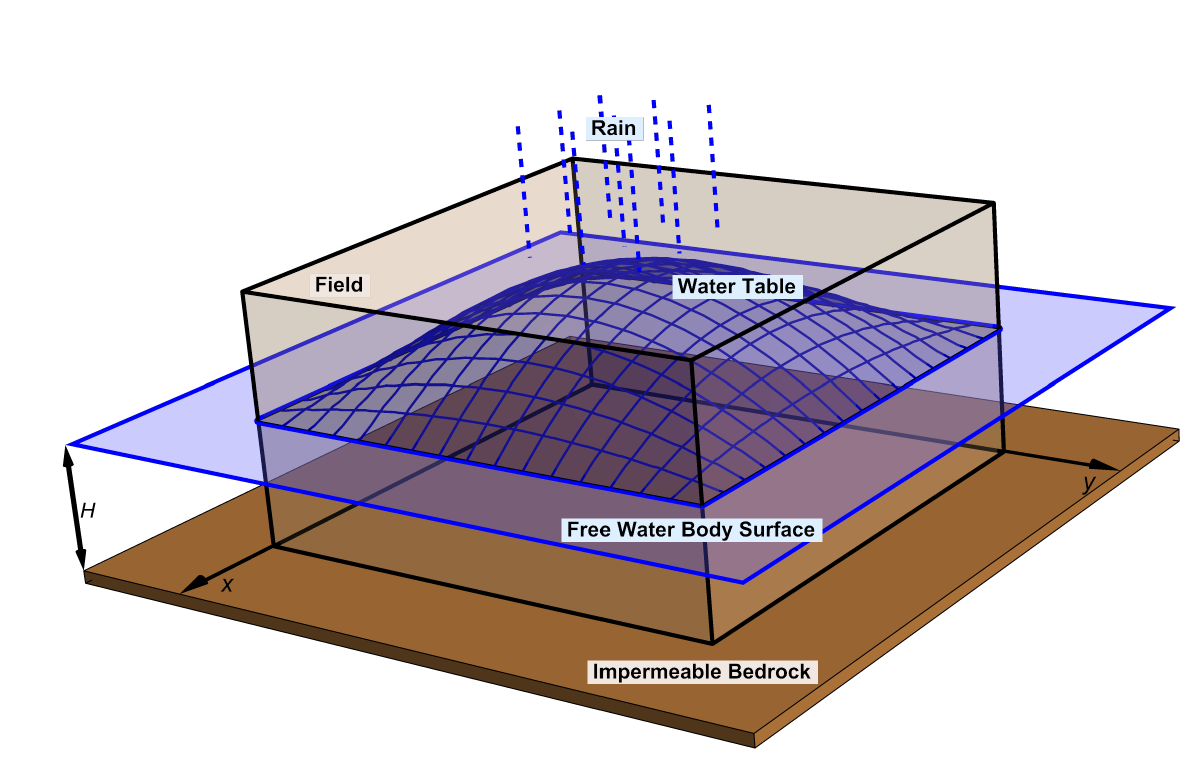}}
\end{picture}
\caption{Fields (brown) drained by a grid of ditches (blue). (a) Topographical view. Domain $\Omega$ represents a field surrounded by a free water body. (b) 3D detail of a field surrounded by free water body accumulated in ditches.}
\end{figure}

\begin{figure}[htp]
\label{fig:parallel}
\setlength{\unitlength}{1cm}
\begin{picture}(12.5,10.5)
\put(0,4.9){
\includegraphics[width=12cm]{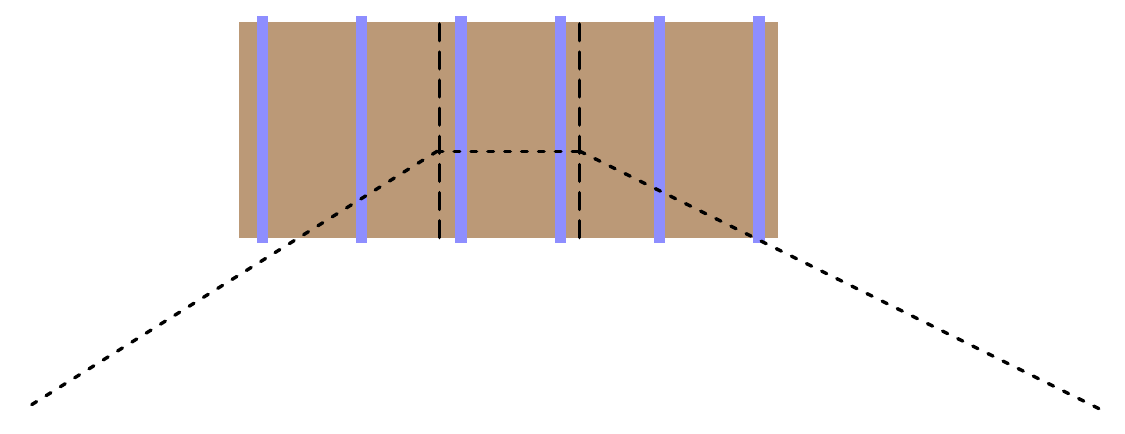}}
\put(4,9.7){(a)}
\put(4,5.8){(b)}
\put(0,0){
\includegraphics[width=12cm]
{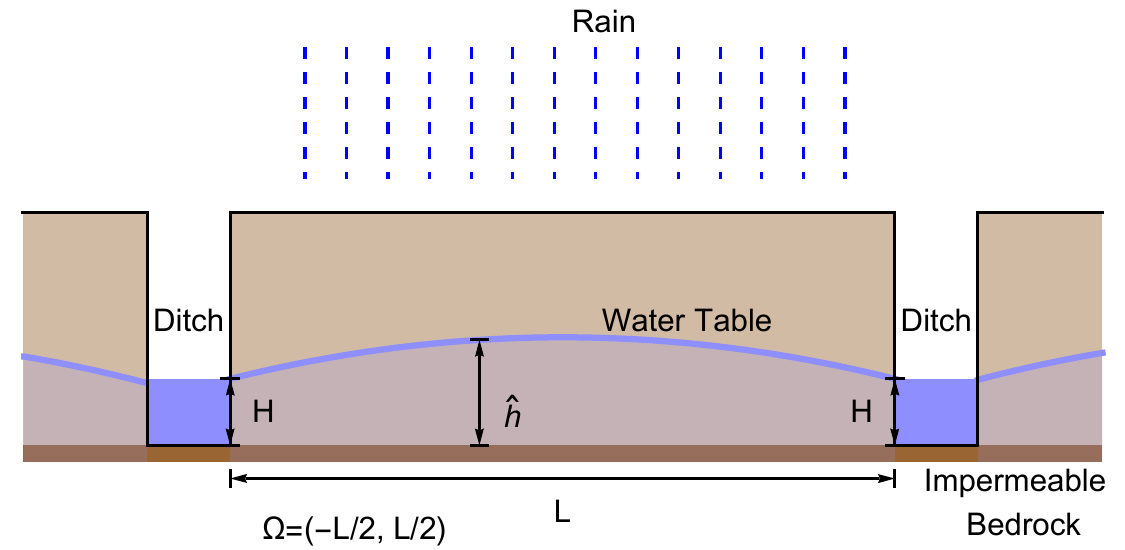}
}
\end{picture}
\caption{Fields (brown) drained by parallel ditches (blue). (a) Topographical view. (b) 2D detail of a field between two ditches.}
\end{figure}

\section{Maximum and Comparison Principles}
\label{sec:maxprincip}
We will address the question of validity of SMP for a nonnegative weak solution to
$$
\left\{\begin{aligned}
\frac{\partial}{\partial t} b(u(\boldsymbol{x}, t))-\Delta_p u(\boldsymbol{x}, t) & =f(\boldsymbol{x}, t) & & \text { for }(\boldsymbol{x}, t) \in \Omega \times(0, T) ; \\
u(\boldsymbol{x}, 0) & =u_0(\boldsymbol{x}) & & \text { for } \boldsymbol{x} \in \Omega ; \\
u(\boldsymbol{x}, t) & =0 & & \text { for }(\boldsymbol{x}, t) \in \partial \Omega \times(0, T) ,
\end{aligned}\right.
$$
which is \eqref{in:prob} with $c = 1$ for simplicity.
The reader is referred, e.g., to \cite{ArielSpringer2021} for the precise definition of weak solution
and the appropriate function spaces.
Recall that
we assume the following hypotheses: $b: \mathbb{R}_{+} \rightarrow \mathbb{R}_{+}$ is a continuous function, $b(0)=0$, and $b \in C^1(0,+\infty)$ with $b^{\prime}>0$ in $(0,+\infty)$. For simplicity, we assume that both, $f: \Omega \times(0, T) \rightarrow \mathbb{R}$ and $u_0: \Omega \rightarrow \mathbb{R}$, are continuous and nonnegative. 
We also assume here that $f \geq 0$ on $\Omega \times(0, T)$. Then indeed, $u \geq 0$ on $\Omega \times(0, T)$ by Weak Comparison Principle (WCP) \cite[Proposition 1, p.~28]{ArielSpringer2021}.

In the modelling of groundwater flow, 
the way the water spreads from the source area (where rain or irrigation occurs) is crucial question. 
In particular, we are interested in the following scenario. 
The initial distribution of water $u_0 \in C\left(\Omega\right)$ is such that $u_0 > 0$ on some compact subset of $\Omega$. 
Considering diffusion process, one would expect that the groundwater immediately spreads toward the boundary $\partial\Omega$. 
Indeed, the solution to linear diffusion equation satisfies the following SMP.
\begin{definition}[Strong Maximum Principle]
    {\rm 
    Let $f \geq 0$ on $\Omega \times (0,T)$ and
    $u \in C\left(\overline{\Omega} \times [0,T)\right)$ 
    be a corresponding nonnegative weak solution to \eqref{in:prob}. 
    We say that $u$ satisfies SMP if there exists $\tau \in (0, T)$ such that  $u(\boldsymbol{x},t) > 0$ for all $(\boldsymbol{x},t) \in \Omega \times (0, \tau)$ and $u(\boldsymbol{x},t) \equiv 0$ for all $(\boldsymbol{x},t) \in \Omega \times [\tau, T)$.
    }
\end{definition}

In the nonlinear case $p \neq 2$ and/or $b(s)\not\equiv s$, the situation is not so clear-cut at all. 
The validity of SMP depends on $p$ and properties of $b(s)$ as $s \to 0+$.
At first, we recall famous Barenblatt's self-similar solution of Leibenson's equation obtained in~\cite{Barenblatt1952b}.

\begin{example}[Barenblatt's self-similar solution]
Let $N\in\mathbb{N}$, $k, p>1$ be given constants.
Existence of self-similar radially symmetric solution of Leibenson's equation
$$
\frac{\partial u^{1/k}}{\partial t}=
\Delta_{p}u\quad 
\mbox{in } \mathbb{R}^N\times(0,+\infty)\,,
$$
was established in \cite{Barenblatt1952b}. In particular, it was shown using ODE techniques that the following radially symmetric problem 
$$
\begin{aligned}
    \frac{\partial}{\partial t} w^{1/k}
    &
    = 
    \frac{1}{r^{N - 1}}\frac{\partial}{\partial r}
    \left[
    r^{N - 1} \left|\frac{\partial w}{\partial r}\right|^{p - 2}
    \frac{\partial w}{\partial r}
    \right]
    & & \text { in } (0,+\infty)\times(0, +\infty)\,, 
\end{aligned}
$$
possesses self-similar solution, i.e., solution of the form
$$
w(r,t)=
t^{-\frac{N}{\beta}}
\mathcal{B}_{N,k,p}
\left(\frac{r}{t^{1/\beta}}\right)\,,
$$
where $\beta=\beta(N,k,p)>0$, which satisfies 
$\int_0^{+\infty} w(r,t) r^{N-1}\mathrm{d}r=1$ for any $t>0$.
Then, for a self-similar radially symmetric solution of Leibenson's equation, it holds $u(\boldsymbol{x},t)= w(|\boldsymbol{x}|,t)$ for $\boldsymbol{x}\in\mathbb{R}^N$ and $t>0$.
Moreover, it easily follows that the initial trace of $u$ is the Dirac measure in $\mathbb{R}^N$concentrated at $0$, i.e.,
$\lim_{t\to 0+} u(\cdot,t)\to \delta_0$
in the sense of measures in $\mathbb{R}^N$.

Three qualitatively different cases need to be distinguished here.
\begin{itemize}
\item[(i)] If $k(p-1)>1$, then $$\mathcal{B}_{N,k,p}(s)=
\left[
\max\left\{0,C-\kappa s^{\frac{p}{p-1}}\right\}\right]^{\gamma}\,,$$ where
$\gamma=\gamma(N, k, p)>0$, $\kappa=\kappa(N, k, p)>0$, $C=C(N, k, p)>0$.
\item[(ii)] If $k(p-1)=1$, then
$\beta=p$ and
$$\mathcal{B}_{N,k,p}(s)=C\,\exp \left(-\left(\frac{s}{p}\right)^{\frac{p}{p-1}}\right)\,,$$
where $C=C(N,k,p)>0$.
\item[(iii)] If $0<k(p-1)<1$, then 
$$\mathcal{B}_{N,k,p}(s)=
\left(C+\kappa s^{\frac{p}{p-1}}\right)^{\gamma}\,,$$ where
$\gamma=\gamma(N, k, p)>0$, $\kappa=\kappa(N, k, p)>0$, $C=C(N, k, p)>0$.
\end{itemize}
It can be easily seen that the self-similar solution $w$ has compact support $\mathop{\mathrm{supp}} w = [0,t^{\frac{1}{\beta }} \left(\frac{C}{\kappa
   }\right)^{\frac{p-1}{p}}]$ in the case (i), while $w(r,t)>0$
   for any $r\geq 0, t>0$, in the cases (ii) and (iii). This makes the case (i) qualitatively distinct from the other two cases. 
   
   Now we turn our eyes back to \eqref{in:prob}. Let
   $k,p>1$,
   $k (p-1) > 1$, $b(s)=s^{1/k}$, and $\boldsymbol{x_0}\in\Omega$ be arbitrary but fixed.
   Let also $0<\sigma<\frac{1}{2}
   {\left(
   \kappa/C\right)^{\beta(p-1)/p}}
   {\mathop{\mathrm{dist}(\boldsymbol{x}_0, \partial\Omega)}}^{\beta}
   $ be fixed.
Then the function 
$$
u(\boldsymbol{x},t)\equiv w(|\boldsymbol{x}-\boldsymbol{x}_0|,t+\sigma)
$$ 
is a solution to \eqref{in:prob} with $f(\boldsymbol{x},t)\equiv 0$ and $u_0(\boldsymbol{x})\equiv w(|\boldsymbol{x}-\boldsymbol{x}_0|,\sigma)$ for all $\boldsymbol{x}\in\Omega$ and 
   $0<t<
   \frac{1}{2}
   \left(
   \kappa/C
   \right)^{\beta(p-1)/p}
   \mathop{\mathrm{dist}(\boldsymbol{x}_0, \partial\Omega)}^{\beta}$. Taking into account that $u_0\geq 0$ and it is positive on a set of positive measure, this means that SMP does not hold for~\eqref{in:prob} for this choice of parameters and function $b$. This motivated our research in \cite{BenediktGirgKotrlaTakac2017, BenediktGirgKotrlaTakac2019} to obtain sufficient conditions for the validity of SMP as well as discovery of further counterexamples for SCP.
\end{example}

Indeed, we obtained the following affirmative result for the SMP, see~\cite{BenediktGirgKotrlaTakac2017} for more details and the proof.

\begin{proposition}
{\rm (see \cite[Thm 1.1.]{BenediktGirgKotrlaTakac2017})} Let $1<p<2, N \geq 1$ and assume that $b: \mathbb{R}_{+} \rightarrow \mathbb{R}_{+}$is as above in \eqref{in:prob} and satisfies also
\begin{equation}
\label{eq:MP:SMP:cond}
\lim _{s \rightarrow 0+} \frac{s^{2-p} b^{\prime}(s)}{|\log s|^{p-1}}=0 .
\end{equation}

Finally, assume that $u: \bar{\Omega} \times[0, T) \rightarrow \mathbb{R}_{+}$ is a continuous, nonnegative, weak solution to \eqref{in:prob}. Then, for any fixed $t_0 \in(0, T)$, the solution $u\left(\cdot, t_0\right)$ is either positive everywhere on $\Omega$ or else identically zero on $\Omega$.

In particular, if $u(\xi, 0)=u_0(\xi)>0$ for some $\xi \in \Omega$, then there exists $\tau \in(0, T]$ such that $u(x, t)>0$ for all $(x, t) \in \Omega \times(0, \tau)$, i.e., Strong Maximum Principle is valid in the $(N+1)$-dimensional space-time cylinder $\Omega \times(0, \tau)$. The number $\tau \in(0, T)$ can be estimated from below by
$$
\tau=\sup \left\{T^{\prime} \in(0, T]: u(\xi, t)>0 \text { for all } t \in\left[0, T^{\prime}\right)\right\}>0
$$
\end{proposition}
Note that $b$ is given by \eqref{eq:bu:water} in our models. 
Hence, $b'(s) > k \equiv {\rm const.}$ for all $s \geq 0$ 
and the validity of SMP depends on value of $p$ only.
The following counterexample is slightly modified \cite[Example 2.3, pp.~368--369]{BenediktGirgKotrlaTakac2019}, where the comparison of stationary solution $u$ and evolutionary solution $v$ to \eqref{in:prob} with $b(s) \equiv s$ is studied in one space dimension. To obtain a counterexample to SMP, we will compare evolutionary solution with the trivial one.

\begin{example}[Counterexample to SMP for $p > 2$ in 1D]
Let $p > 2$, $\max\{2,p/(p -2)\} < \beta$ and $0 < \gamma \leq \beta$. Assume also that $b'(s) \geq k \equiv {\rm const.} > 0$. 
In accordance with \cite{BenediktGirgKotrlaTakac2019}, define
\begin{alignat}{2}
\nonumber
u(x,t)&\equiv u(x) = 0\, \qquad \text{ and } \\
v(x,t)&= 
\left\{\hspace{2.00mm}
\begin{aligned}
& 0 && \quad\mbox{for }x \in (-1,0]\,, \\
&t\, x^{\beta}\,(1-x^{\gamma}) && \quad\mbox{for }x \in (0,1)\,.
\end{aligned}
\right.
\end{alignat}
Clearly, $u(\pm 1,t) = v(\pm 1, t) = 0$, $v(x,0) = u(x,0) = 0$, $u(x,t) = v(x,t) = 0$ on $(-1,0) \times (0,T)$, and $0 = u(x,t) < v(x,t)$ on $(0,1) \times (0,T)$.

Our goal is to show that there exists $t_0 > 0$ such that
$$
f(x,t) 
\stackrel{\rm def}{=}b'(v)\frac{\partial v}{\partial t} 
- 
\frac{\partial}{\partial x}\left(\left|\frac{\partial v}{\partial x}\right|^{p - 2} \frac{\partial v}{\partial x}\right) > 0 \quad \mbox{ on } (0,1) \times (0,t_0)\,.
$$
Indeed,
$$
\begin{aligned}
f(x,t) 
& 
\geq k \frac{\partial v}{\partial t} 
- 
\frac{\partial}{\partial x}\left(\left|\frac{\partial v}{\partial x}\right|^{p - 2} \frac{\partial v}{\partial x}\right) 
\\
&
=
kx^{\beta}
\left[
(1 - x^{\gamma}) - (p - 1)\beta^{p - 1}(\beta - 1) k^{-1} t^{p - 1} x^{(\beta - 1)(p - 2) - 2}
\times 
\right.
\\
&
\times
\left.
\left|
1 - \frac{\beta + \gamma}{\beta} x^{\gamma}
\right|^{p - 2} 
\left(
1 - \frac{(\beta + \gamma)(\beta + \gamma - 1)}{\beta(\beta - 1)} x^{\gamma}
\right)
\right]\,.
\end{aligned}
$$
If 
$$
     \left(\frac{\beta(\beta - 1)}{(\beta + \gamma)(\beta + \gamma - 1)}\right)^{1/\gamma}
     < 
     x
     <
     1\,,
$$
then
$$
    \left(
    1 - \frac{(\beta + \gamma)(\beta + \gamma - 1)}{\beta(\beta - 1)} x^{\gamma}
    \right)
    < 0
$$
which together with 
$(1 - x^{\gamma}) > 0$ ensures  that $f(x,t) > 0$.
On the other hand, if 
$$
0 < x
\leq
\left(\frac{\beta(\beta - 1)}{(\beta + \gamma)(\beta + \gamma - 1)}\right)^{1/\gamma}\,,
$$
we have $(1 - x^{\gamma}) > k_1 \equiv {\rm const.}$ and, hence, we may find $t_0$ small enough such that 
$$
\begin{aligned}
&
\left[ 
(1 - x^{\gamma}) - (p - 1)\beta^{p - 1}(\beta - 1) k^{-1} t^{p - 1} x^{(\beta - 1)(p - 2) - 2}
\times 
\right.
\\
&
\times
\left.
\left|
1 - \frac{\beta + \gamma}{\beta} x^{\gamma}
\right|^{p - 2} 
\left(
1 - \frac{(\beta + \gamma)(\beta + \gamma - 1)}{\beta(\beta - 1)} x^{\gamma}
\right)
\right] > 0\,.
\end{aligned}
$$
Recall that the constants $\beta$ and $\gamma$ are chosen such that $(\beta - 1)(p - 2) - 2 > 0$.
This concludes the counterexample to SMP.
\end{example}

The counterexample above provides the existence of a nonnegative weak solution 
with compact support $[0,1]$ for sufficiently small times. Let us observe that the solution in our counterexample is qualitatively different from the Barenblatt's solution, since its compact support stays constant for sufficiently small times. Moreover, positive bump of our solution is induced by suitable nontrivial nonnegative right-hand side of the equation rather than by the initial condition, which is identically zero in our counterexample.

In \cite{BenediktGirgKotrlaTakac2019}, we studied counterexamples to SCP with weak solutions $u$ and $v$ that are positive in $(-1,1)$ 
such that the only point satisfying $u'(x) = v'(x) = 0$ being $x = 0$.
Thus $x = 0$ was the only point of degeneracy of the diffusive part driven by the $p$-Laplacian with $p>2$. Such degeneracy causes that a perturbation from one side of the point of degeneracy $x=0$ does not spread through this point immediately in general, but some waiting time is needed. This result from \cite{BenediktGirgKotrlaTakac2019} is adopted to \eqref{in:prob}
in the following proposition.

\begin{proposition}
	Let $2 < p < \infty$, $\max\{2, p/(p - 2)\} \leq \alpha < \beta$, $\gamma = \beta - \alpha$ $({}>0)$, 
    $b'(s) > k \equiv {\rm const.}$ for all $s \geq 0$, and set 
	\begin{align*}
		 u(x,t) \equiv u(x) &= 1 - |x|^{\alpha} \qquad \mbox{ together with }\\
		 f_s(x) \stackrel{\rm def}{=} - \frac{\partial}{\partial x}\left(\left|\frac{\partial u}{\partial x}\right|^{p - 2} \frac{\partial u}{\partial x} \right) & = (p - 1) \alpha^{p - 1}(\alpha - 1) |x|^{(\alpha - 1)(p - 2) + \alpha - 2}\,,  
	\end{align*}
 both for $x \in (-1,1)$.
	Furthermore, let
	$$
		h_i(x,t) = f(x) + |x|^{\beta} \psi_i(x,t)  \quad \mbox{ for }	(x,t) \in (-1,1) \times (0, T)\,; \ i=1,2,	
	$$
	where $\psi_i \in L^{\infty} \left( (-1,1) \times (0,T) \right)$ satisfy 
	\begin{itemize}
		\item[(i)] $0 \leq \psi_1 \leq \psi_2 \leq \frac{1}{2}|x|^{\beta}(1 - |x|^{\gamma} + C t)$ with some constant $C > 0$; and
		\item[(ii)] $\psi_1 \not\equiv \psi_2$ on $(-1,1)\times(0, \tau)$ for any $\tau \in (0,T)$. 
	\end{itemize}
	Let $w_i(x,t)$ $(i = 1,2)$ be the weak solutions to 
	\begin{equation}
	\label{eq:w}
	\left\{\hspace{2.00mm}
	\begin{alignedat}{2}
			\textstyle\frac{\partial b(w_i)}{\partial t} - \frac{\partial}{\partial x}\left(\left|\frac{\partial w_i}{\partial x}\right|^{p-2}\frac{\partial w_i}{\partial x}\right) & = h_i(x,t) &&\quad\mbox{ for } (x,t)\in (-1,1)\times (0, T)\,; \\
			w_i(\pm 1, t)&= 0 &&\quad\mbox{ for } t\in (0, T)\,; \\
			w_i(x,0)&=u_0(x) &&\quad\mbox{ for } x\in (-1,1)  
	\end{alignedat}
	\right.
	\end{equation}	
	which satisfy $w_i \in C\left([-1,1] \times [0, T]\right)$. 
	If the constant $C > 0$ is sufficiently small, there exists $t_0 > 0$ such that 
	
	\begin{equation}
	\label{eq:conclusion}
	\left\{ \hspace*{2.00mm}		
	\begin{aligned}
		w_1(x,t) & \leq w_2(x,t) && \qquad \mbox{ for } (x,t) \in (-1,1) \times (0, t_0)\,,\\
		w_1 & \not\equiv w_2&& \qquad \mbox{ on } (-1,1) \times (0, t_0)\,,\\
		\mbox{but } w_1(0,t) &= w_2(0,t) && \qquad \mbox{ for {\em all} } t \in  (0,t_0)\,.
	\end{aligned}
	\right.	
	\end{equation}
	
	Finally, if $\psi_1 \equiv \psi_2 \equiv 0$ on $(-1, 0) \times (0, t_0)$, then we have also $w_1 \equiv w_2 \equiv u$ on $(-1,0] \times (0, t_0)$ 
	and $w_1 \not\equiv w_2$ on $(0,1) \times (0, t_0)$.
\end{proposition} 

The proof follows the same steps as \cite[Proof of Theorem 2.4, pp.~369--370]{BenediktGirgKotrlaTakac2019}. 
Define
\begin{equation*}
\hat v(x,t) =u(x) + t\, |x|^{\beta}\,(1-|x|^{\gamma}) \quad\mbox{for }(x,t) \in (-1,1)\times(0,T)
\end{equation*}
and
\begin{equation*}
\hat g(x,t) = b'(\hat{v})\frac{\partial \hat{v}}{\partial t} - \frac{\partial}{\partial x}\left(\left|\frac{\partial \hat{v}}{\partial x}\right|^{p - 2} \frac{\partial \hat{v}}{\partial x} \right) \,. 
\end{equation*}
Since $b'(s) \geq k$ for all $s \geq 0$,
$$
\begin{aligned}
\hat{g}(x,t) - f_s(x) 
&
= 
b'(\hat{v})\frac{\partial \hat{v}}{\partial t} - \frac{\partial}{\partial x}\left(\left|\frac{\partial \hat{v}}{\partial x}\right|^{p - 2} \frac{\partial \hat{v}}{\partial x} \right) + \frac{\partial}{\partial x}\left(\left|\frac{\partial u}{\partial x}\right|^{p - 2} \frac{\partial u}{\partial x} \right)\\
&
\geq
k\frac{\partial \hat{v}}{\partial t} - \frac{\partial}{\partial x}\left(\left|\frac{\partial \hat{v}}{\partial x}\right|^{p - 2} \frac{\partial \hat{v}}{\partial x} \right) + \frac{\partial}{\partial x}\left(\left|\frac{\partial u}{\partial x}\right|^{p - 2} \frac{\partial u}{\partial x} \right)
\end{aligned}
$$
and we obtain by similar calculation as in 
\cite[Example 2.3, pp.~368--369]{BenediktGirgKotrlaTakac2019}
that there exists $t_0 \in (0,T)$ such that 
$$
    \hat g(x,t) - f_s(x) \geq \frac{1}{2}|x|^{\beta}\left(1 - |x|^{\gamma} + \frac{C}{k}t\right)
    \quad \mbox{ for } (x,t) \in (-1,1) \times (0, t_0)\,,
$$  
provided $C > 0$ is chosen small enough. We refer reader to \cite[Proposition 1, p.~28]{ArielSpringer2021} for the Weak Comparison Principle (WCP) applicable to problems with $b(s) \not\equiv s$. This ends the outline of the proof.

Results presented in this section demonstrate that the validity of SMP for quasilinear parabolic problems is a complex issue. 
The relevance of SMP gains importance in connection with the fact that explicit forms of solutions to problem~\eqref{in:prob} can be found only in very exceptional cases and we are limited to the use of numerical methods.
While SMP might not hold in general, classical methods such as implicit Euler or Crank-Nicholson methods remain reasonable effective when it is known to be valid. 
However, for solutions with compact support, specialized numerical methods are needed to mitigate unwanted numerical dispersion at the boundary of the compact supports and avoid oscillations of the numerical solution, see, e.g.,~\cite{
ARBOGAST2019108921, Berger1979, Cavalli20072098,
GU2020109378, Jiang2021, Liu2011939, Parlange2000339}.
Our further research on the validity of SMP aims to provide criteria for selecting accurate and efficient numerical methods for a given problem based on SMP validity or presence of compact support solutions.

\section{Parameter estimation and CFD simulations of experiments}
\label{sec:CFD}

Understanding and managing groundwater flow in fractured hard-rock aquifers are crucial for a variety of purposes, 
including sustainable water access in rural areas and dewatering construction projects like tunnels, see, e.g., \cite{Gustafson1994, larsson1984ground, Shapiro2002}.
However, accurately estimating parameters in constitutive laws for these aquifers presents significant challenges compared to traditional porous media.
These parameters are usually obtained by fitting formulas such as, e.g., \eqref{const:Darcy} or \eqref{power:law} on data obtained by series of laboratory experiments performed on samples of given porous medium, see, e.g., \cite{Bear1972, Darcy1856, King1898}.
This process is suitable for porous media encompassing unconsolidated materials like soils, sands, and gravels, or permeable rocks like sandstones.
In the context of fractured hard rock, however, this approach encounters significant limitations. Hard rocks have extremely low permeability and the groundwater flow occurs almost exclusively within the network of fractures. It may occur that the network of fractures can have a low density of significant fractures, making it difficult to obtain representative samples for traditional laboratory experiments. 
Such experiments would require large, sometimes even multi-cubic meter rock samples to ensure enough fractures are captured, necessitating careful extraction techniques to minimize damage to the natural fracture network. This poses significant logistical and cost limitations. 
For fracture systems, a common approach is to conduct physical experiments or numerical simulations on a single fracture or a small number of intersecting fractures, allowing the influence of fracture intersections to manifest within the constitutive relationships. This strategy has been adopted in numerous studies across experimental, theoretical, and numerical domains, see, e.g., \cite{Berkowitz2002, Brush2003, cherubini2013evidence, LapcevicNovakowskiSudicky, 
QUINN2020124384,
sarkar2004fluid, StarkVolker,  yan2018non, Zimmerman2004163}.

\begin{figure}[htb]
	\centering	\includegraphics[width=8cm]{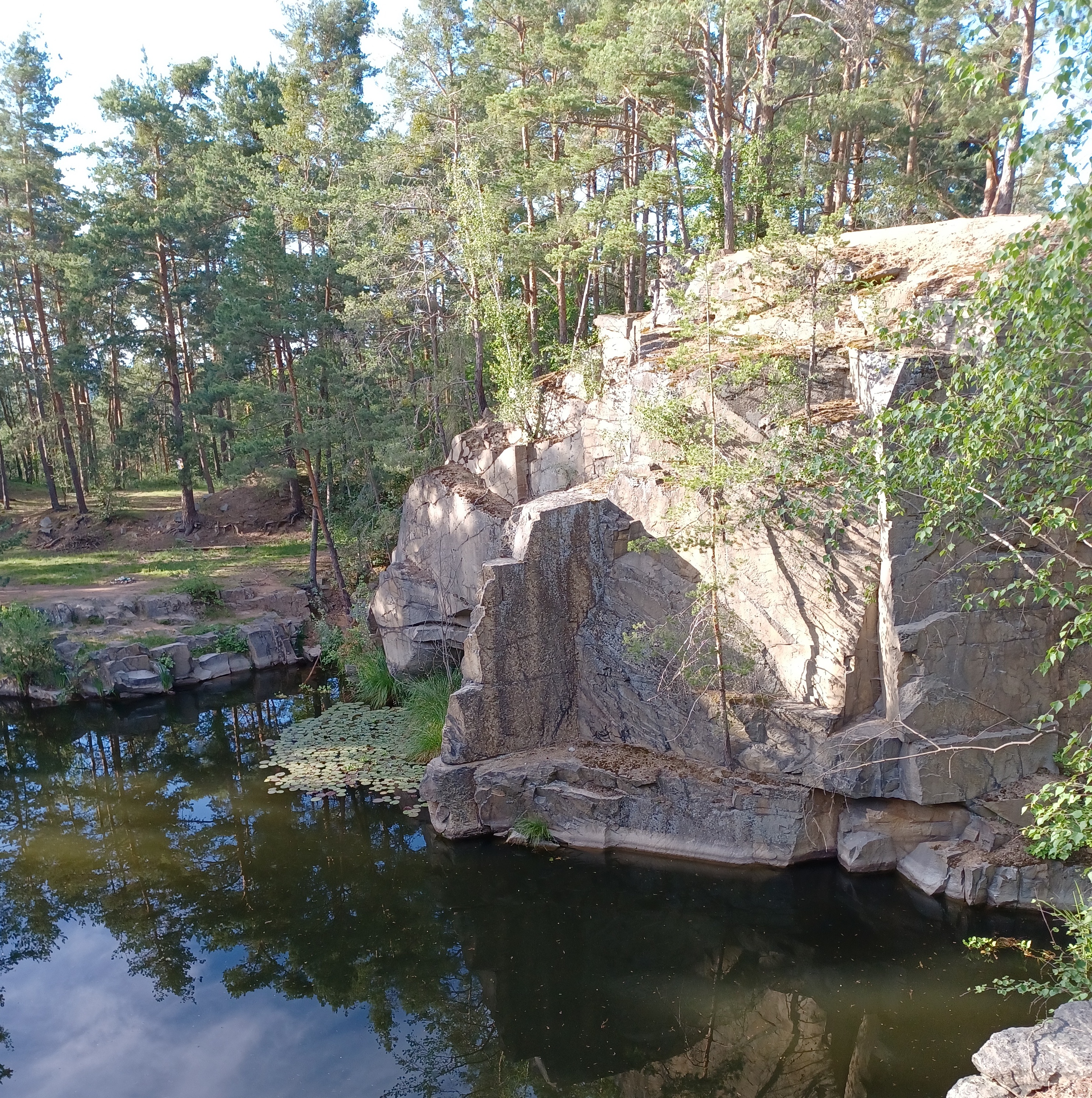}\\
	\caption{Granite rock exposures in abandoned quarry Špic by Něčín, Czech Republic.}
	\label{fig: lomSpic}
\end{figure}

While laboratory experiments remain valuable, approach based on computer fluid dynamics (CFD) simulations of physical experiments offers a potentially cheaper and more accessible alternative. By simulating the flow in several intersecting fractures, we aim to estimate the parameters.
To obtain as realistic outputs as possible, we created 2D and 3D geometrical models, see Fig.~\ref{fig: 3Dsharp_2D_pukliny} and~\ref{fig: 3Dsharp_3D_pukliny}, of several intersecting fractures found on easy accessible granite rock exposures in abandoned and partially flooded quarry Špic by Něčín, Czech Republic, see Fig.~\ref{fig: lomSpic}. 

\begin{figure}[htb]
	\centering
\includegraphics[width=12cm]{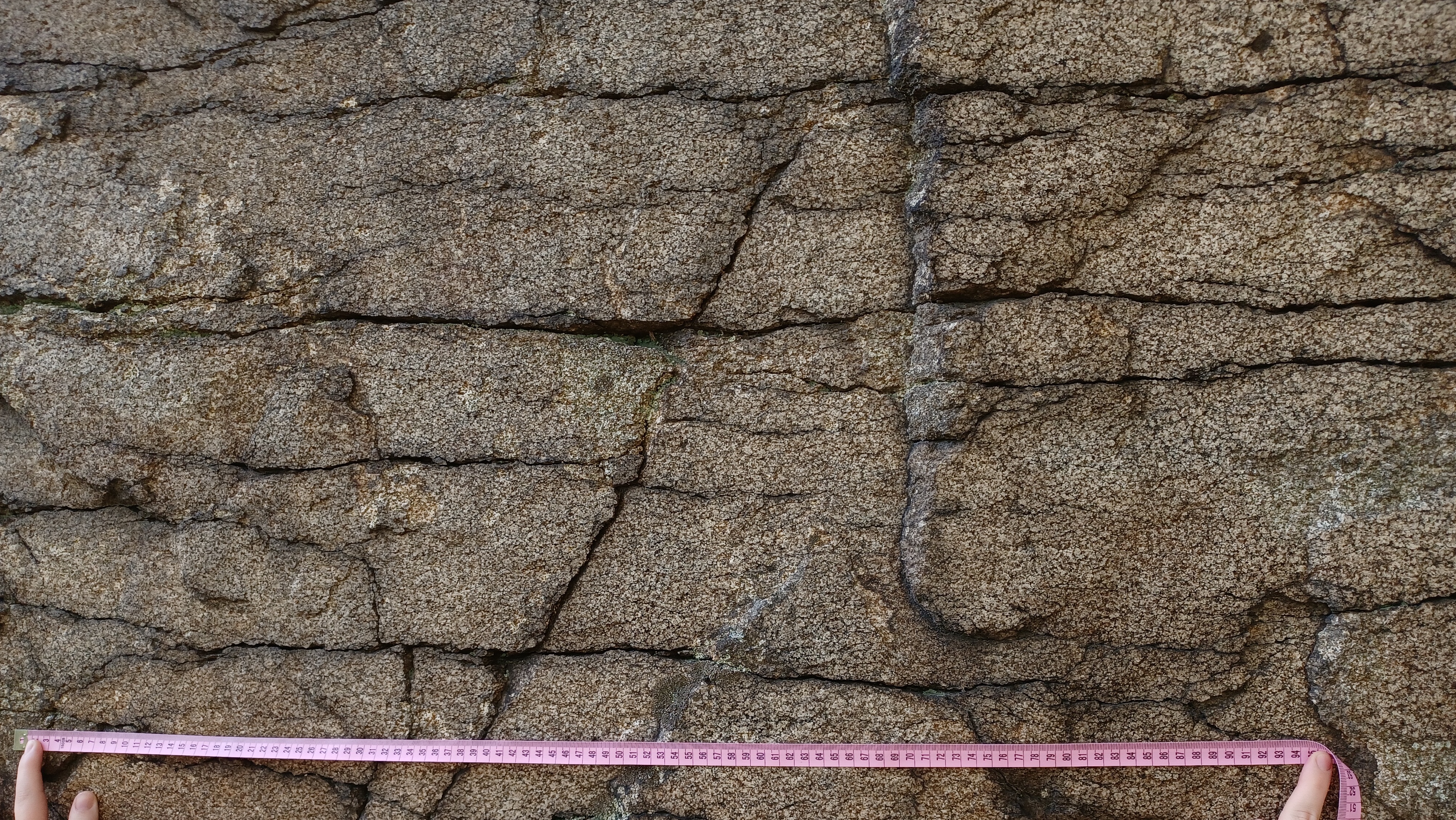}\\
	\caption{Fracture network in the granite rock exposure, quarry Špic by Něčín.}
	\label{fig: lomSpic_pukliny}
\end{figure}

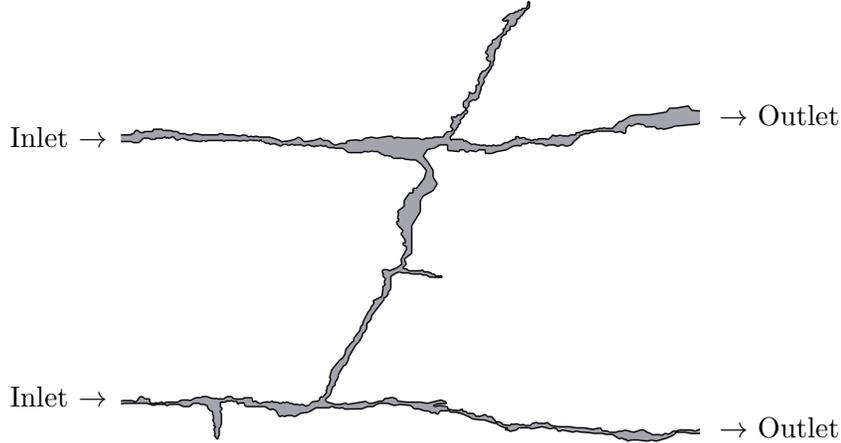
\begin{figure}[htb]
\setlength{\unitlength}{1cm}
\begin{picture}(11.5,7.5)
\centering
\includegraphics[width=10cm]{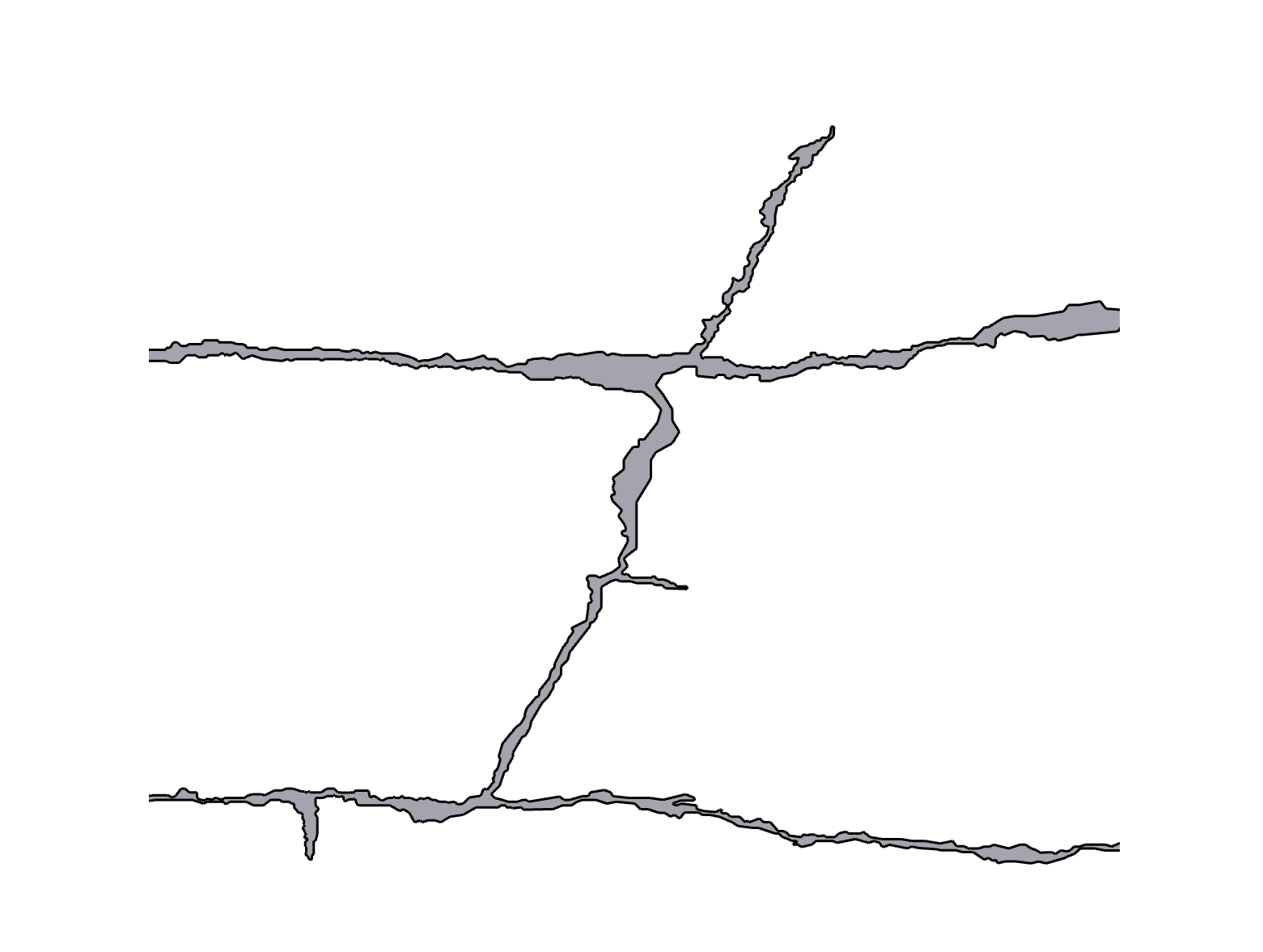}
\put(-10.3,1.15){Inlet $\rightarrow$}
\put(-10.3,4.61){Inlet $\rightarrow$}
\put(-1,0.75){ $\rightarrow$ Outlet}
\put(-1,4.9){ $\rightarrow$ Outlet}
\end{picture}
	\caption{2D model of fracture network based on the most significant fractures from the photograph on Fig.~\ref{fig: lomSpic_pukliny}.}
	\label{fig: 3Dsharp_2D_pukliny}
\end{figure}

\begin{figure}[htb]
	\centering
	\includegraphics[width=5.7cm]{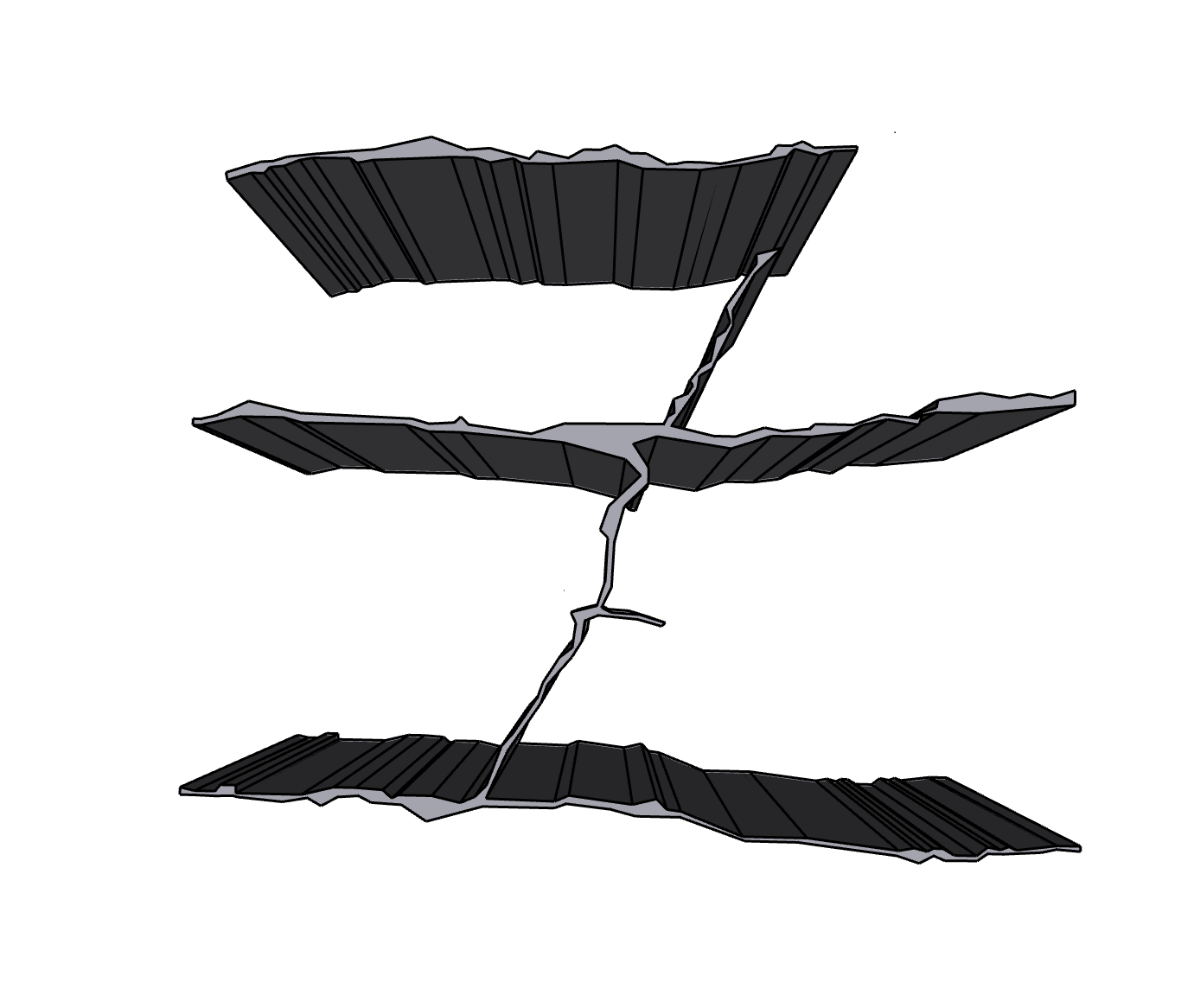}
	\includegraphics[width=5.7cm]{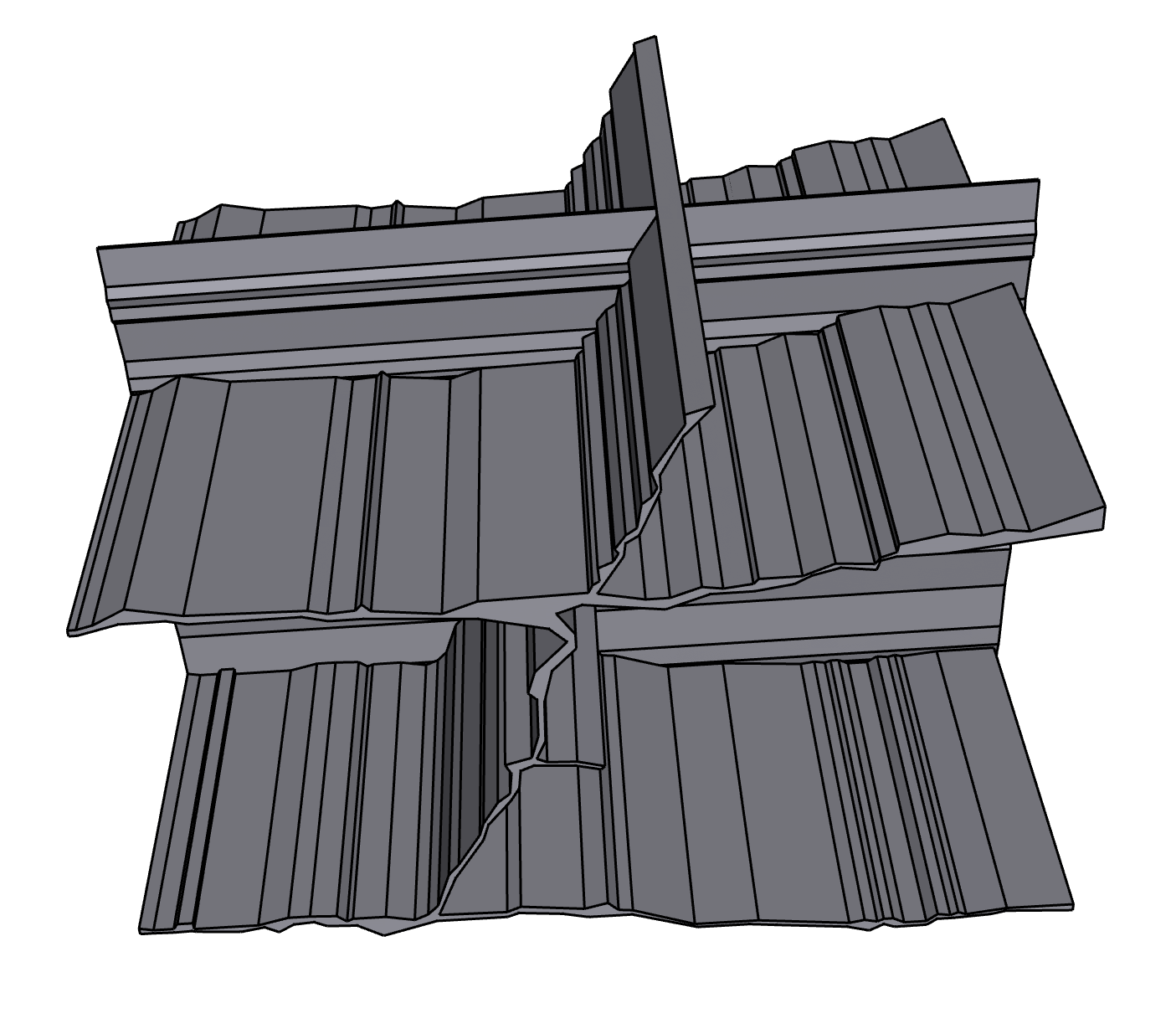}
	\caption{3D model 
 of fracture network. (a) Fracture network based on the most significant fractures from the photograph on Fig.~\ref{fig: lomSpic_pukliny}. (b) 3D model with artificially added vertical fracture.}
	\label{fig: 3Dsharp_3D_pukliny}
\end{figure}

Using these geometrical models, see \cite{diplomka}, we simulated  physical
experiments to obtain datasets suitable for parameter estimation by numerical solving the stationary incompressible Navier-Stokes equation together with the continuity equation
\begin{samepage}
\begin{eqnarray}
\label{eq:NS}
-\nu \Delta\vec{v}+
\vec{v} \cdot \nabla \vec{v} &=&-\nabla P \\ 
\label{eq:cont}
\nabla \cdot \vec{v}&=&0
\end{eqnarray}
\end{samepage}
in the fracture network represented by the domain $\Omega_{F2D} \subset\mathbb{R}^2$
(Fig.~\ref{fig: 3Dsharp_2D_pukliny})
or $\Omega_{F3D} \subset \mathbb{R}^3$
(Fig.~\ref{fig: 3Dsharp_3D_pukliny}).
Here $\vec{v}$
is the velocity field ($\Omega_{F2D}\to\mathbb{R}^2$ or 
$\Omega_{F3D}\to\mathbb{R}^3$), $P$ is the pressure field, and $\nu$ is the kinematic viscosity. The system \eqref{eq:NS}-\eqref{eq:cont} is completed by imposing boundary conditions, see \eqref{eq:bc:2D} and \eqref{eq:bc:3D} for 2D and 3D case, respectively.

Let us note that our aim was to test feasibility of this approach. 
We do not claim that obtained results from this preliminary stage are ready to be used in practice, see discussion in Remark~\ref{rem:anisotropy}. 

\subsection{2D model of fracture network}

Our 2D model of fracture network is based on the most significant fractures from the photograph on Fig.~\ref{fig: lomSpic_pukliny} and is represented by domain  $\Omega_{\mathrm{F2D}}\subset(0, 0.55)\times(0, 0.5)$ (in meters) with Lipschitz boundary $\partial\Omega_{\mathrm{F2D}}$. 
This boundary
$\partial\Omega_{\mathrm{F2D}}=
\Gamma_{\mathrm{inlet}}\cup
\Gamma_{\mathrm{outlet}}\cup
\Gamma_{\mathrm{wall}}$
is a union of pairwise disjoint sets, where $\Gamma_{\mathrm{inlet}}$, $\Gamma_{\mathrm{outlet}}$, and
$\Gamma_{\mathrm{wall}}$ represent inlet into the fracture network, outlet from the fracture network, and fixed walls of the fracture network, respectively.
Our 2D model of fracture network has two inlets and two outlets to take into account possible mixing effects inside the network, see Fig.~\ref{fig: 3Dsharp_2D_pukliny}. 
More formally, we have   
$\Gamma_{\mathrm{inlet}}=
\{(x,z)\in
\partial\Omega_{\mathrm{F2D}}
\colon x=x_{\mathrm{inlet}},
z\in I_{\mathrm{in}_1}
\cup I_{\mathrm{in}_2}
\}$ 
and 
$\Gamma_{\mathrm{outlet}}=\{
(x,z)\in
\partial\Omega_{\mathrm{F2D}}
\colon x=x_{\mathrm{outlet}},
z\in I_{\mathrm{out}_1}
\cup I_{\mathrm{out}_2}
\}$, where $x_{\mathrm{inlet}}=0$,
$x_{\mathrm{outlet}}=0.55$,
$I_{\mathrm{in}_j},
I_{\mathrm{out}_j} \subset (0, 0.5)$, $j=1,2$,
are open intervals and
$
I_{\mathrm{in}_1}\cap I_{\mathrm{in}_2}=\emptyset$, 
$I_{\mathrm{out}_1}\cap I_{\mathrm{out}_2}=\emptyset$.
Thus, outer normal vector fields $\vec{n}$ on $\Gamma_{\mathrm{inlet}}$ and
$\Gamma_{\mathrm{outlet}}$
are constant fields, 
$\vec{n}=(-1,0)$ and 
$\vec{n}=(1,0)$, respectively.

The flow in the domain $\Omega_{\mathrm{F2D}}$ was simulated in \emph{OpenFOAM} using \emph{simpleFOAM} solver by numerically solving \eqref{eq:NS}-\eqref{eq:cont} with boundary conditions
\begin{equation}
\label{eq:bc:2D}
\left\{
	\begin{array}{ll}
    \displaystyle
	\vec{v}=(0,0)\quad 
    &\mbox{on } \Gamma_{\mathrm{wall}}\,,\\[0.2cm]
\displaystyle    
\vec{v}=(v_{\mathrm{inlet}},0)\quad &\mbox{on }\Gamma_{\mathrm{inlet}}\,,
\\[0.2cm]
\displaystyle
P=0,\quad v\mbox{ satisfies condition described below}
& \mbox{on }\Gamma_{\mathrm{outlet}}\,.
\end{array}
\right.
\end{equation}	
We impose \emph{fluxCorrectedVelocity}
outflow condition
provided by \emph{OpenFOAM}
on $\Gamma_{\mathrm{outlet}}$.
In essence, this nonlocal condition acts computationally in the following way. An initial estimate of the velocity at the outlet is obtained so that the change in velocity across the outlet boundary is zero. This initial guess is then corrected based on the calculated flux leaving the domain. This correction ensures that the flux leaving the domain matches the expected flow based on the pressure and the internal flow field. For detailed information, see the \emph{OpenFOAM} documentation~\cite{OF}.

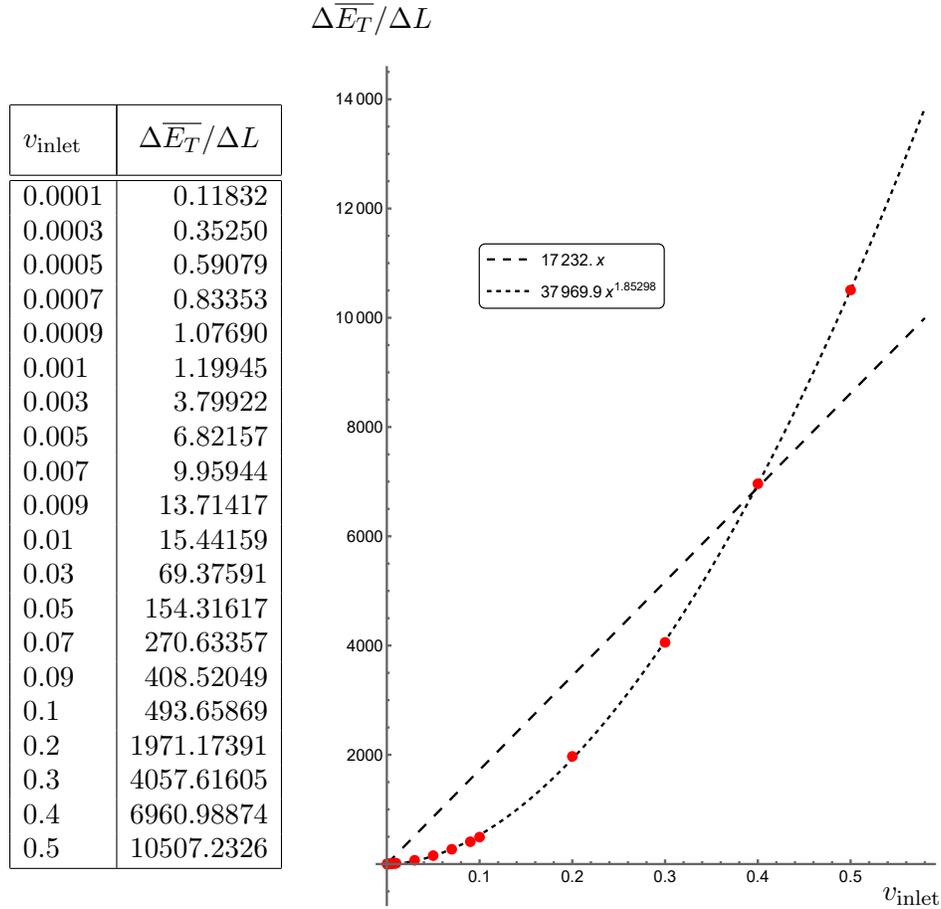
\begin{figure}[htp]
\setlength{\unitlength}{1cm}
\begin{picture}(12.5,12.5)
\put(4.2,0){
\includegraphics[width=8cm]{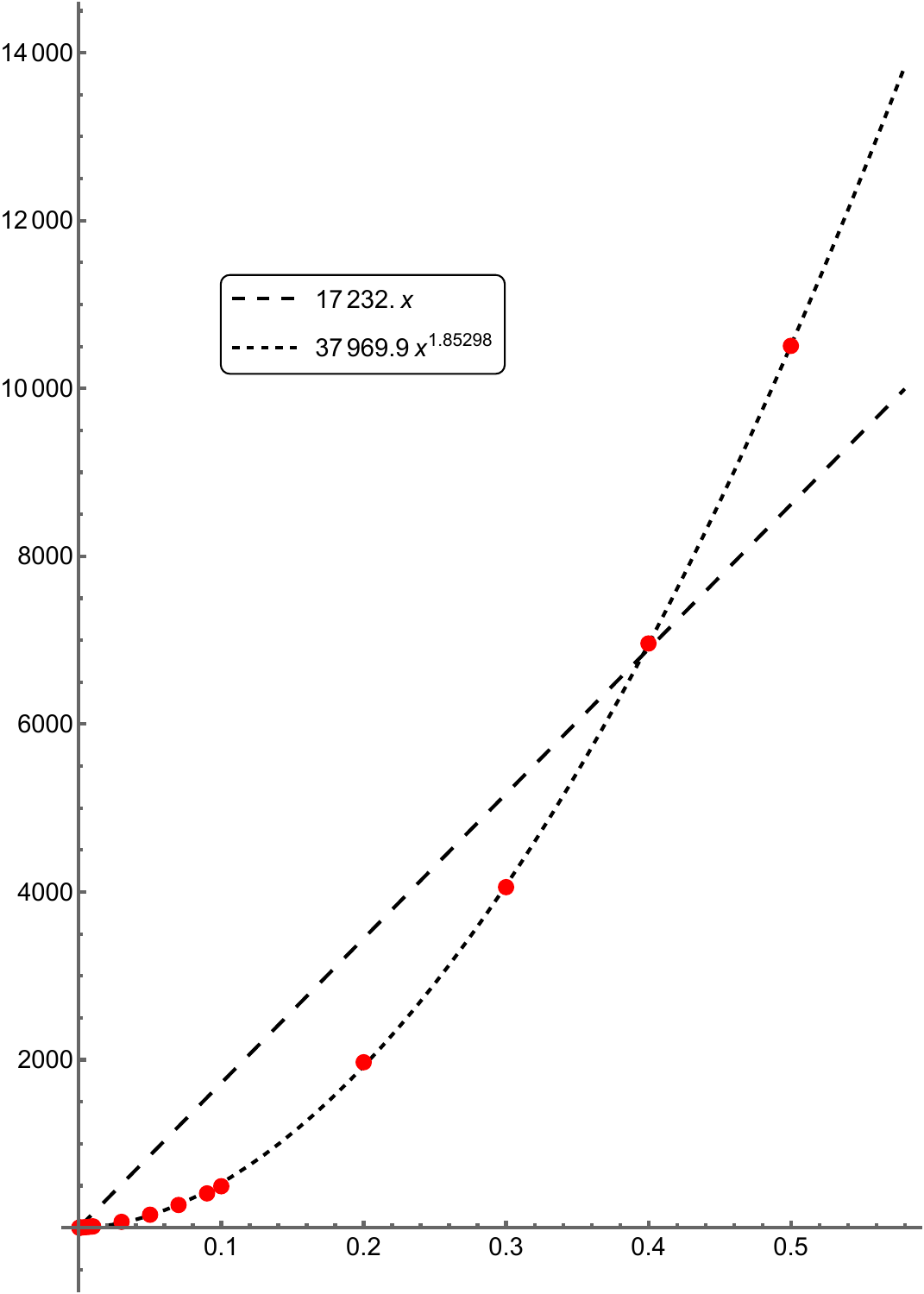}}
\put(4,11.7){$\Delta \overline{E_T} / \Delta L$}
\put(11.6,0.1){$v_{\mathrm{inlet}}$}
\put(0,5.5){\begin{tabular}{|l|r|}  
\hline
$v_{\mathrm{inlet}}$ & \multicolumn{1}{|c|}{$\vphantom{\biggl(}\Delta \overline{E_T} / \Delta L$} \\
\hline 
\hline 
0.0001 & 0.11832 \\
0.0003 & 0.35250 \\
0.0005 & 0.59079 \\
0.0007 & 0.83353 \\
0.0009 & 1.07690 \\
0.001  & 1.19945 \\
0.003  & 3.79922 \\
0.005  & 6.82157 \\
0.007  & 9.95944 \\
0.009  & 13.71417 \\
0.01   & 15.44159 \\
0.03   & 69.37591 \\
0.05   & 154.31617 \\
0.07   & 270.63357 \\
0.09   & 408.52049 \\
0.1    & 493.65869 \\
0.2    & 1971.17391 \\
0.3    & 4057.61605 \\
0.4    & 6960.98874 \\
0.5    & 10507.2326 \\
\hline
\end{tabular}}
\end{picture}
\caption{Table: collected data from numerical simulations performed on 2D model. Graphs: Darcy's and power type law fitted to the data.}
\label{fig:fitted}
\end{figure}

Simulation is performed for several values of  $v_{\mathrm{inlet}}$, see Figure~\ref{fig:fitted}. 
In principle, this is a numerical imitation of Darcy's physical experiment with an adaptation for a fractured rock. In his physical experiment, Darcy determined the total mechanical energy loss 
$E_T$ using the hydraulic head $h$, which neglects kinetic energy but can be easily measured in reality using piezometers. In a numerical experiment, we do not have the possibility to measure the hydraulic head using piezometers, so we have to choose a different approach. Just as in a physical experiment, we want to determine the mechanical energy loss of the fluid flowing through the fracture network. In our numerical experiment, we are modeling the actual flow in the fractures and we work with actual velocity $\vec{v}$ not with the averaged velocity $\vec{v}_{\mathrm{av}}$. In our case, the term representing kinetic energy $1/2\,\rho\,v^2$ is not negligible, where $v$ stands for the magnitude of $\vec{v}$. 

Since the velocity and pressure are dependent on location, we need to consider average total mechanical energy per unit volume. For this purpose, we introduce  
\begin{eqnarray}
\displaystyle
\overline{E_T}_{\mathrm{inlet}} &=& 
\frac{
\int_{\widehat{\Gamma}_{\mathrm{inlet}}}
E_T\mathop{
\mathrm{d}z}
}{
\int_{{\widehat{\Gamma}_{\mathrm{inlet}}}}
\mathop{
\mathrm{d}z}}\,, \\
\overline{E_T}_{\mathrm{outlet}} &=&
\frac{
\int_{\widehat{\Gamma}_{\mathrm{outlet}}}
E_T\mathop{
\mathrm{d}z}
}{
\int_{{\widehat{\Gamma}_{\mathrm{outlet}}}}
\mathop{
\mathrm{d}z}}\,,
\end{eqnarray}
where
\begin{eqnarray*}
\widehat{\Gamma}_{\mathrm{inlet}}&=&
\{z\in\mathbb{R}\colon (x_{\mathrm{inlet}},z)\in
\Gamma_{\mathrm{inlet}}\} = I_{\mathrm{in}_1}
\cup I_{\mathrm{in}_2}\,,\\
\widehat{\Gamma}_{\mathrm{outlet}}&=&
\{z\in\mathbb{R}\colon (x_{\mathrm{outlet}},z)\in
\Gamma_{\mathrm{outlet}}\} 
=
I_{\mathrm{out}_1}
\cup I_{\mathrm{out}_2}\,.
\end{eqnarray*}
Similar averaging approach was used in \cite[p.~297]{LUCAS2007295} to define macroscopic pressure and velocity.

Now, we obtain the constitutive relationship from numerical simulations as follows. We test the model for a range of  values $v_\mathrm{inlet}$ in~\eqref{eq:bc:2D}
and record the values of expression 
$$\frac{\triangle\overline{E_T}}{\triangle L}=\frac{\overline{E_T}_{\mathrm{inlet}}-\overline{E_T}_{\mathrm{outlet}}}{\triangle L}\,.$$
Then we approximate 
dependence 
$$
\frac{\triangle\overline{E_T}}{\triangle L} = f_{\mathrm{law}}(u_{\mathrm{inlet}})
$$
by fitting parameters in the expression for $f_{\mathrm{law}}$  to simulated flow data.
We used the following expressions 
\begin{eqnarray}   
\label{f:Darcy}
f_{\mathrm{law}}(u_{\mathrm{inlet}})&=& \alpha\, u_{\mathrm{inlet}}\,,\mbox{ with parameter }\alpha>0,\mbox{ for Darcy's law}\,,\\[0.2cm]
\label{f:power}
f_{\mathrm{law}}(u_{\mathrm{inlet}})&=& \beta\, u_{\mathrm{inlet}}^{\gamma}\,,
\mbox{ with parameters }
\beta,\gamma>0,\mbox{ for power law}.  
\end{eqnarray}
The constitutive law for specific discharge $q=A_{\mathrm{inlet}} u_{\mathrm{inlet}}$ is then 
$$
q= A_{\mathrm{inlet}} f_{\mathrm{law}}^{-1}\left(\frac{
\triangle \overline{E_T}}{\triangle L}\right)\,,
$$
where 
$A_{\mathrm{inlet}}= \int_{\widehat{\Gamma}_{\mathrm{inlet}}}\mathrm{d}z$ is the cross-section area of the inlet.

Now we will make a~connection of results of these simulations to groundwater flow in fractured hard-rock aquifers.  The water table is observed in vertical boreholes. Due to the much larger diameter of the borehole compared to the fractures, the groundwater flow experiences a sudden expansion upon entering. This rapid increase in cross-sectional area significantly reduces flow velocity and so the contribution of the term corresponding to kinetic energy to total mechanical energy per volume in the borehole can be neglected.
Thus the water would rise approximately to level $h_{\mathrm{inlet}}=
\overline{E_T}_{\mathrm{inlet}}/
(\varrho g)$ and $h_{\mathrm{outlet}}=\overline{E_T}_{\mathrm{outlet}}/(\varrho g)$ above the vertical datum if the fictive boreholes are located in the inlet and outlet of the fracture network, respectively.
This leads us to the following form of the constitutive law
\begin{equation}
\label{konst:H:T}
q = A_{\mathrm{inlet}} f_{\mathrm{law}}^{-1}
\left(\varrho\, g\,
\frac{
h_{\mathrm{inlet}}-h_{\mathrm{outlet}}
}{\triangle L}\right)\,.
\end{equation}

By fitting each law~\eqref{f:Darcy}
and~\eqref{f:power}
to data collected from numerical simulations, we obtained
$\alpha=17232.004$ and $\beta=37969.871$, $\gamma=1.852982$, with significantly lower root mean square error $17.1324$ for the power-type law~\eqref{f:power} compared to $776.396$ for Darcy's law~\eqref{f:Darcy}, see Fig.~\ref{fig:fitted}.
This finding strongly suggests that the power-type law provides a much more accurate description of fluid flow behavior in the specific type of fracture networks simulated in this study. 
Using these fitted values of $\beta$ and $\gamma$ in~\eqref{konst:H:T} together with $A_{\mathrm{inlet}}=9.9987\cdot 10^{-3}$ (computed for our geometrical model), and $g=9.8066$ (conventional standard value for gravitational acceleration)
\begin{equation}
\label{const:our:fracture2D}
q = 0.004815
\left(
\frac{
h_{\mathrm{inlet}}-h_{\mathrm{outlet}}
}{\triangle L}
\right)^{0.5397}
\,.
\end{equation}

Assuming that the fractured hard-rock aquifer is homogeneous, isotropic (in the sense of porosity due to fractures) and formed by a system of fractures of the same type as in the studied section, 
we can use the two-dimensional differential form \eqref{vec:power:law} of constitutive law \eqref{const:our:fracture2D} for substitution in the balance equation \eqref{eq:balance}
to obtain the equation for the groundwater level
\begin{equation}
\label{eq:doubly:hath:2Dsim}
0.0347\,\frac{\partial \widehat{h}}{\partial t}-
0.004815\,\operatorname{div}
\left(\widehat{h} |\nabla \widehat{h}|^{-0.4603}\nabla \widehat{h}\right)=\widehat{g}(x, y, t)\,,
\end{equation}
where we have used calculated value $\phi_{\mathrm{eff}}=0.0347$ (the ratio of volume the fracture network and volume of the sample of the rock) for our 2D geometric model of fracture network.
Let us note that the term $|\nabla \widehat{h}|^{-0.4603}\nabla \widehat{h}$
is understood in the sense that it returns zero vector for zero vector as input, cf~\eqref{vec:power:law}.

\begin{remark}
\label{rem:anisotropy}
Let us note that the hard-rock aquifers are often nor isotropic due to prevailing orientation of the fractures  nor homogeneous due to varying density and aperture of the fractures, see, e.g., \cite{KIRALY1971255, larsson1984ground, Marechal2004}. In our future research, we plan to address these topics, especially the issue of anisotropy of nonlinear constitutive laws. 
\end{remark}

\subsection{3D model of fracture network.}
Our 3D fracture network model also incorporates the most significant fractures from the photograph on Fig.~\ref{fig: lomSpic_pukliny}.  To achieve a spatially representative network, a single artificial fracture was added and connected to existing fractures within the rock mass. Although the precise location was not based on a specific observation, it aligns with typical fracture distributions observed in the granitic exposures of the quarry.
The fracture model is then represented by domain  $\Omega_{\mathrm{F3D}}\subset (0,0.55)\times (0, 0.5)\times (0, 0.5)$ with Lipschitz boundary $\partial\Omega_{\mathrm{F3D}}$. This boundary $\partial\Omega_{\mathrm{F3D}}=
\Gamma_{\mathrm{inlet}}\cup
\Gamma_{\mathrm{outlet}}\cup
\Gamma_{\mathrm{wall}}$
is a union of pairwise disjoint sets, where $\Gamma_{\mathrm{inlet}}$, $\Gamma_{\mathrm{outlet}}$, and
$\Gamma_{\mathrm{wall}}$ represent inlet into the fracture network, outlet from the fracture network, and fixed walls of the fracture network, respectively.
Our 3D model of fracture network, see Fig.~\ref{fig: 3Dsharp_3D_pukliny}, is such that
$\Gamma_{\mathrm{inlet}}\subset\{
(x,y,z)\in
\partial\Omega_{\mathrm{F3D}}
\colon x=x_{\mathrm{inlet}}
\}$ 
and 
$\Gamma_{\mathrm{outlet}}\subset\{
(x,y,z)\in
\partial\Omega_{\mathrm{F3D}}
\colon x=x_{\mathrm{outlet}}
\}$ are sufficiently regular, where $x_{\mathrm{inlet}}=0$
and $x_{\mathrm{outlet}}=0.55$, so that outer normal vector fields $\vec{n}$ are well defined on $\Gamma_{\mathrm{inlet}}$ and
$\Gamma_{\mathrm{outlet}}$, respectively.
Moreover, they
are constant fields, 
$\vec{n}=(-1,0,0)$ and 
$\vec{n}=(1,0,0)$
on $\Gamma_{\mathrm{inlet}}$ and
$
\Gamma_{\mathrm{outlet}}$, respectively. 

\begin{figure}[htp]
\setlength{\unitlength}{1cm}
\begin{picture}(12.5,13.2)
\put(4.2,0){
\includegraphics[width=8cm]{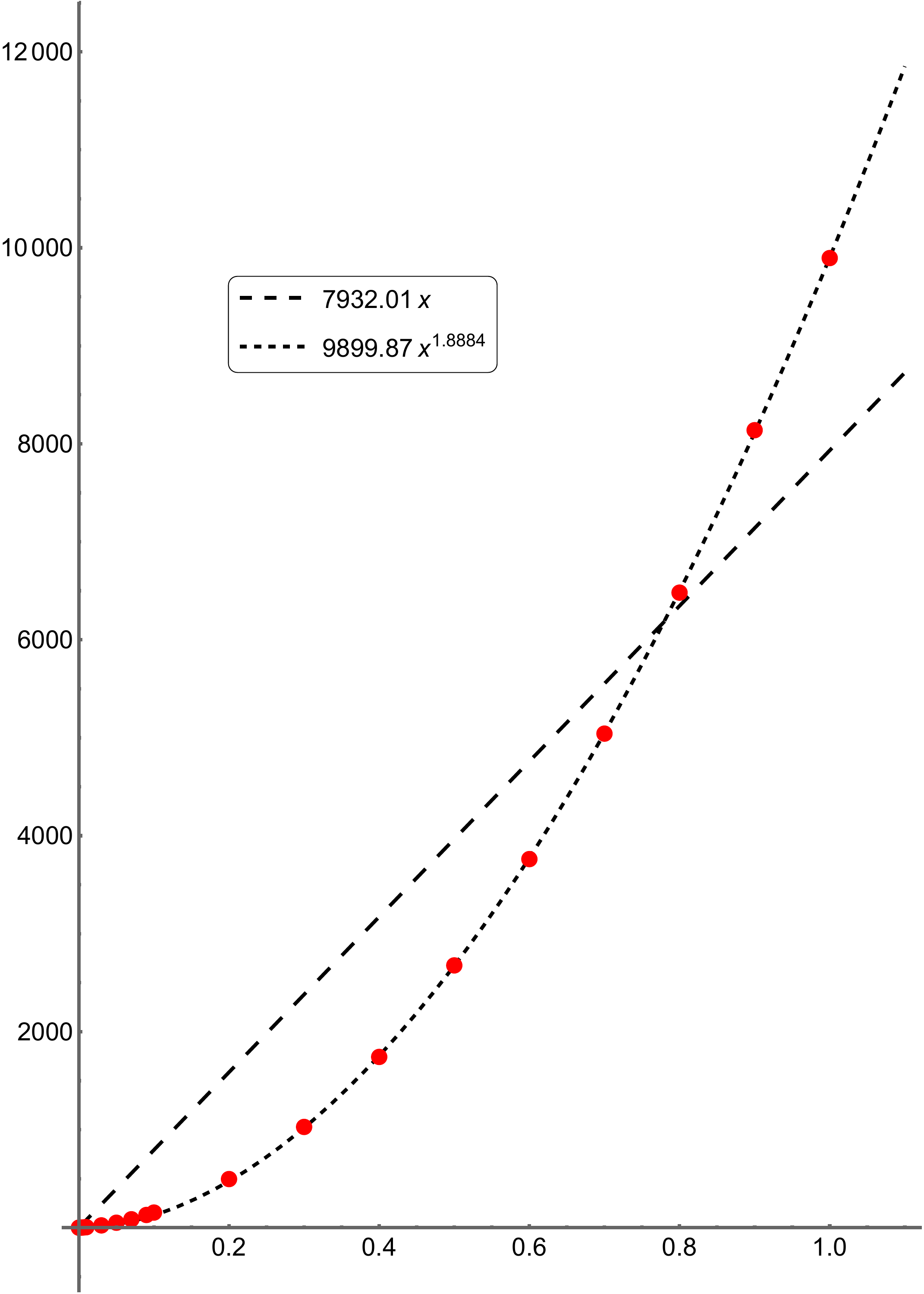}}
\put(4,11.7){$\Delta \overline{E_T} / \Delta L$}
\put(11.6,0.1){$v_{\mathrm{inlet}}$}
\put(0,6.15){\begin{tabular}{|l|r|}  
\hline
$v_{\mathrm{inlet}}$ & \multicolumn{1}{|c|}{$\vphantom{\biggl(}\Delta \overline{E_T} / \Delta L$} \\
\hline 
\hline 
0.0001 &  0.03358\\
0.0003 &  0.10272\\
0.0005 &  0.17364\\
0.0007 &  0.2462\\
0.0009 &  0.3228\\
0.001 &   0.36216\\
0.003 &   1.20726\\
0.005 &   2.21733\\
0.007 &   3.3564\\
0.009 &   4.60464\\
0.01 &    5.26071\\
0.03 &    23.40721\\
0.05 &    50.52472\\
0.07 &    85.79307\\
0.09 &   129.60561\\
0.1 &    154.31271\\
0.2 &    496.05733\\
0.3 &   1028.10647\\
0.4 &   1742.67221\\
0.5 &   2676.70004\\
0.6 &   3761.46226\\
0.7 &   5041.91289\\
0.8 &   6480.71039\\
0.9 &   8138.64571\\
1.0 &   9895.35284\\
\hline
\end{tabular}}
\end{picture}
\caption{Table: collected data from numerical simulations performed on 3D model. Graphs: Darcy's and power type law fitted to the data.}
\label{fig:fitted3D}
\end{figure}
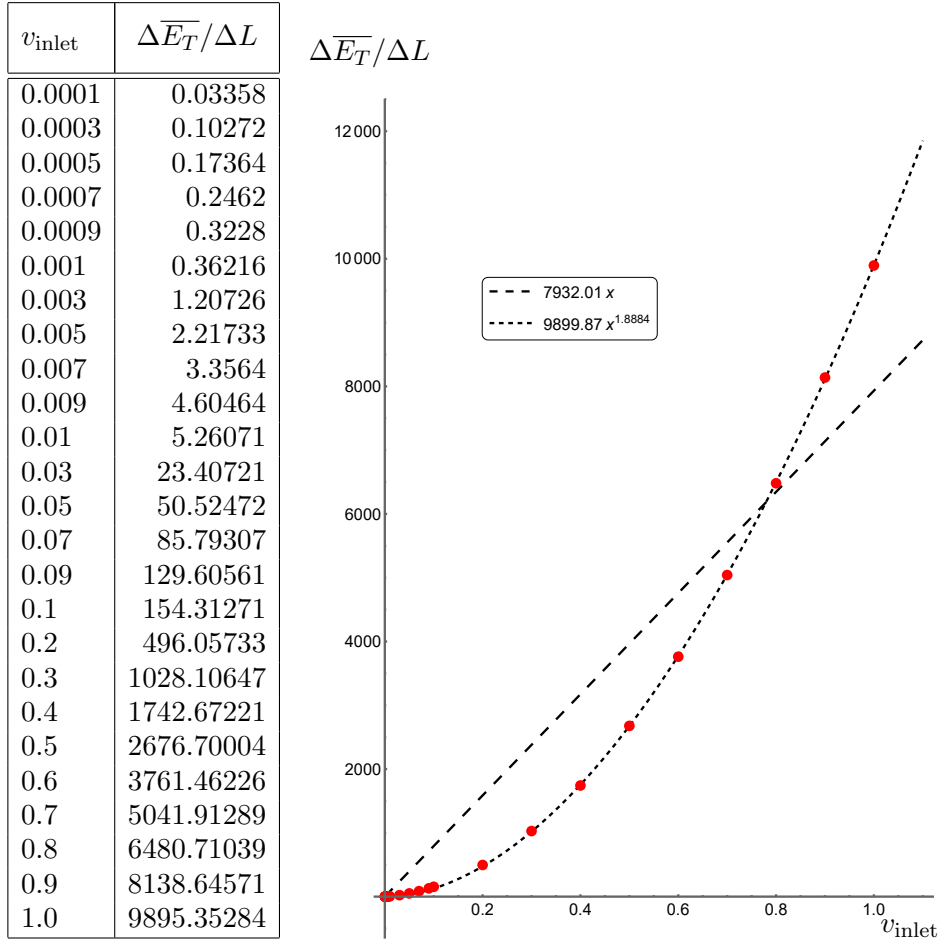

The flow in the domain $\Omega_{\mathrm{F3D}}$ was simulated in \emph{OpenFOAM} using \emph{simpleFOAM} solver by numerically solving \eqref{eq:NS}-\eqref{eq:cont} with boundary conditions
\begin{equation}
\label{eq:bc:3D}
\left\{
	\begin{array}{ll}
    \displaystyle
	\vec{v}=(0,0,0)\quad 
    &\mbox{on } \Gamma_{\mathrm{wall}}\,,\\[0.2cm]
\displaystyle    
\vec{v}=(v_{\mathrm{inlet}},0,0)\quad &\mbox{on }\Gamma_{\mathrm{inlet}}\,,\\[0.2cm]
\displaystyle
P=0,\quad \mbox{ \emph{fluxCorrectedVelocity} condition for } \vec{v} &\mbox{on }\Gamma_{\mathrm{outlet}}\,.
\end{array}
\right.
\end{equation}	
Simulation was performed for several values of  $v_{\mathrm{inlet}}$, see Figure~\ref{fig:fitted3D}. The procedure was analogous as in the 2D case, using averaged values of total energy over the surfaces of inlet and outlet.
By fitting each law~\eqref{f:Darcy}
and~\eqref{f:power}, we obtained
$\alpha=7932.011$ and 
$\beta=9899.873$, $\gamma=1.888399$.
Again, we found that root mean square error $12.44$ for the power-type law~\eqref{f:power} 
is significantly lower compared to $748.39$ for Darcy's law~\eqref{f:Darcy}, see Fig.~\ref{fig:fitted3D}. This strongly suggests that the power-type law is much more accurate
for the specific type of fracture networks simulated in this study.

Using these fitted values of $\beta$ and $\gamma$ in~\eqref{konst:H:T} together with $A_{\mathrm{inlet}}=0.006$ (computed for our 3D geometrical model), and $g$ as above, we find 
\begin{equation}
\label{const:our:fracture2D:3Dsim}
q = 0.00597
\left(
\frac{
h_{\mathrm{inlet}}-h_{\mathrm{outlet}}
}{\triangle L}
\right)^{0.529549}
\,.
\end{equation}
Under the same assumption on the 
hard-rock aquifer (homogeneous, isotropic, and formed by a system of fractures of the same type as in the studied section), we obtain the equation for the groundwater level
\begin{equation}
\label{eq:doubly:hath:3dsim}
0.04728\,\frac{\partial \widehat{h}}{\partial t}-
0.00597\,\operatorname{div}
\left(\widehat{h} |\nabla \widehat{h}|^{-0.4704}\nabla \widehat{h}\right)=\widehat{g}(x, y, t)\,,
\end{equation}
where we have used calculated value $\phi_{\mathrm{eff}}=0.04728$ (the ratio of the volume of the fracture network and the volume of the sample of the rock) for our 3D geometric model of fracture network.

\section{P. Dr\'abek and the $p$-Laplacian}
\label{sec:drabek}

Prior to embarking on the main discourse of this section, all the three authors of this paper would like to 
extend their sincere congratulations to their esteemed teacher, Pavel 
Dr\'{a}bek, on the occasion
of his 70th birthday.
They would like to express their profound gratitude for his guidance
and introduction to nonlinear analysis, with a particular emphasis on the
$p$-Laplacian, since the early stages of their studies and scientific careers. One of the key themes of P.~Drábek
is the quest for an analogue of the Fredholm alternative for nonlinear operators, particularly for the $p$-Laplacian. 
This topic was a major focus of research at Prague School of Nonlinear Analysis\footnote{Prague School of Nonlinear Analysis (or Prague School for short) was an informal research group of mathematicians primarily from Charles University Prague, the Czechoslovak Academy of Sciences, and other academic institutions based in Prague."} in the 1970s, where it was probably brought by
J.~Nečas, Ph.D. advisor of S.~Fu\v{c}\'{i}k 
(who was later mentor and advisor of P.~Dr\'{a}bek). The groundbreaking results of the Prague School were published in monograph {\it Spectral analysis of nonlinear
operators} by S.~Fučík, J.~Nečas, J.~Souček, V.~Souček, published in 1973, see \cite{FNSS}.
After the premature death of S.~Fu\v{c}\'{\i}k in 1979, this topic gradually lost its importance in the Prague School. Among other things, because it was a very
difficult topic and after the publication of the above mentioned monograph no
more significant breakthroughs were achieved.
Fortunately, this topic did not completely disappear thanks to
P. Drábek. His passion for this topic was ignited under the mentorship of S.~Fu\v{c}\'{i}k, and he has been diligently researching it ever since his diploma thesis. After accomplishing his C.Sc. degree (equivalent of Ph.D.) at Czechoslovak Academy of Sciences, P.~Dr\'{a}bek relocated to Plze\v{n}, where he got position in the department of mathematics at the College of Mechanical and Electrical Engineering (precursor of today's University of West Bohemia). Here he continued to pursue his research in Fredholm alternative for nonlinear operators and the $p$-Laplacian while the Prague School's primary research focus shifted to Navier-Stokes equations at the prompting of 
J.~Ne\v{c}as. P.~Dr\'{a}bek always returned to the topic with a certain time lag and still does. During almost half a century, an interesting and extensive series of articles has been written, from which
each representing a major advance.
In addition, he has long been attracted to this subject
attention, and so some of the seminal articles
have been written without his direct input (co-authorship), but it would hardly be without his persistence in presenting papers at conferences, seminars and in discussions with a number of eminent mathematicians.

\section*{Acknowledgements}
P.~Girg, L. Kotrla and A.~\v{S}vandov\'{a} were supported by the~Grant Agency of~the~Czech
Republic, Grant No.~22-18261S.
Authors acknowledge the helpful assistance of Gemini (formerly Bard), a large language model from Google AI, in enhancing the language quality of this paper.

\def\cprime{$'$}


\begin{thebibliography}{99}

\bibitem{AravinNumerov}
{\sc Aravin, V.~I., and Numerov, S.~N.}
\newblock {\em ``Teoriya dvizheniya zhidkostei i gazov v nedeformiruemoi poristoi srede''}.
\newblock Gosudarstv. Izdat. Tehn.-Teor. Lit., Moscow, 1953.
\newblock English Transl. by A.~Moscona: {\sl ``Theory of Fluid Flow in Undeformable Porous Media''}, Israel Program for Scientific Translations, Jerusalem, 1965.

\bibitem{ARBOGAST2019108921}
{\sc Arbogast, T., Huang, C.-S., and Zhao, X.}
\newblock Finite volume {WENO} schemes for nonlinear parabolic problems with degenerate diffusion on non-uniform meshes.
\newblock {\em Journal of Computational Physics 399\/} (2019), 108921.

\bibitem{AYRAUD20082686}
{\sc Ayraud, V., Aquilina, L., Labasque, T., Pauwels, H., Molenat, J., Pierson-Wickmann, A.-C., Durand, V., Bour, O., Tarits, C., {Le Corre}, P., Fourre, E., Merot, P., and Davy, P.}
\newblock Compartmentalization of physical and chemical properties in hard-rock aquifers deduced from chemical and groundwater age analyses.
\newblock {\em Applied Geochemistry 23}, 9 (2008), 2686--2707.

\bibitem{Barenblatt1952b}
{\sc Barenblatt, G.~I.}
\newblock On self-similar motions of compressible fluid in a porous medium.
\newblock {\em Akad. Nauk SSSR. Prikl. Mat. Meh. 16}, 6 (1952), 679--698.
\newblock In Russian.

\bibitem{Barenblatt1952}
{\sc Barenblatt, G.~I.}
\newblock On some unsteady motions of a liquid and gas in a porous medium.
\newblock {\em Akad. Nauk SSSR. Prikl. Mat. Meh. 16\/} (1952), 67--78.
\newblock In Russian.

\bibitem{BASAK1979}
{\sc Basak, P.}
\newblock An analytical solution for the transient ditch drainage problem.
\newblock {\em Journal of Hydrology 41}, 3 (1979), 377--382.

\bibitem{Bear1972}
{\sc Bear, J.}
\newblock {\em Dynamics of Fluids in Porous Media}.
\newblock Enviromental science series. American Elsevier Publishing Company, Inc., New York, 1972.

\bibitem{Bear2014}
{\sc Bear, J.}
\newblock {\em Dynamics of Fluids in Porous Media}.
\newblock Dover Civil and Mechanical Engineering Series. Dover Publications, Inc., New York, 2014.

\bibitem{ArielSpringer2021}
{\sc Benedikt, J., Girg, P., and Kotrla, L.}
\newblock Nonlinear models of the fluid flow in porous media and their methods of study.
\newblock In {\em Functional differential equations and applications, FDEA-2019. Proceedings of the 7th international conference, Ariel, Israel, September 22--27, 2019}. Singapore: Springer, 2021, pp.~15--42.

\bibitem{BenediktGirgKotrlaTakac2017}
{\sc Benedikt, J., Girg, P., Kotrla, L., and Tak{\'a}{\v{c}}, P.}
\newblock The strong maximum principle in parabolic problems with the {$p$}-{L}aplacian in a domain.
\newblock {\em Appl. Math. Lett. 63\/} (2017), 95--101.

\bibitem{BenediktGirgKotrlaTakac2018}
{\sc Benedikt, J., Girg, P., Kotrla, L., and Tak\'{a}\v{c}, P.}
\newblock Origin of the {$p$}-{L}aplacian and {A}. {M}issbach.
\newblock {\em Electron. J. Differential Equations\/} (2018), Paper No. 16, 17 pp.

\bibitem{BenediktGirgKotrlaTakac2019}
{\sc Benedikt, J., Girg, P., Kotrla, L., and Tak\'{a}\v{c}, P.}
\newblock The strong comparison principle in parabolic problems with the {$p$}-{L}aplacian in a domain.
\newblock {\em Appl. Math. Lett. 98\/} (2019), 365--373.

\bibitem{Berger1979}
{\sc Berger, Alan~E., B. H. R. J. C.~W.}
\newblock A numerical method for solving the problem $u_t - \delta f (u) = 0$.
\newblock {\em ESAIM: Mathematical Modelling and Numerical Analysis - Modélisation Mathématique et Analyse Numérique 13}, 4 (1979), 297--312.

\bibitem{Berkowitz2002}
{\sc Berkowitz, B.}
\newblock Characterizing flow and transport in fractured geological media: A review.
\newblock {\em Advances in water resources 25}, 8-12 (2002), 861--884.

\bibitem{Brush2003}
{\sc Brush, D.~J., and Thomson, N.~R.}
\newblock Fluid flow in synthetic rough-walled fractures: Navier-stokes, stokes, and local cubic law simulations.
\newblock {\em Water Resources Research 39}, 4 (2003).

\bibitem{Cavalli20072098}
{\sc Cavalli, F., Naldi, G., Puppo, G., and Semplice, M.}
\newblock High-order relaxation schemes for nonlinear degenerate diffusion problems.
\newblock {\em SIAM Journal on Numerical Analysis 45}, 5 (2007), 2098 – 2119.

\bibitem{CHENG2024}
{\sc Cheng, H., Wang, F., Li, S., Guan, X., Yang, G., Cheng, Z., Yu, C., Yuan, Y., and Feng, G.}
\newblock Effect of movability of water on the low-velocity pre-darcy flow in clay soil.
\newblock {\em Journal of Rock Mechanics and Geotechnical Engineering\/} (2024).

\bibitem{cherubini2013evidence}
{\sc Cherubini, C., Giasi, C., and Pastore, N.}
\newblock Evidence of non-darcy flow and non-fickian transport in fractured media at laboratory scale.
\newblock {\em Hydrology and Earth System Sciences 17}, 7 (2013), 2599--2611.

\bibitem{Darcy1856}
{\sc Darcy, H.}
\newblock {\em Les fontaines publiques de la ville de Dijon}.
\newblock Victor Dalmont, Paris, 1856.

\bibitem{Forchheimer1901}
{\sc Forchheimer, P.}
\newblock Wasserbewegung durch {B}oden.
\newblock {\em Zeit. Ver. Deutsch. Ing. 45\/} (1901), 1736--1741 and 1781--1788.

\bibitem{FNSS}
{\sc Fu\v{c}\'{\i}k, S., Ne\v{c}as, J., Sou\v{c}ek, J., and Sou\v{c}ek, V.}
\newblock {\em Spectral analysis of nonlinear operators}.
\newblock Lecture Notes in Mathematics, Vol. 346. Springer-Verlag, Berlin-New York, 1973.

\bibitem{Australia2014}
{\sc {Geoscience Australia}}.
\newblock \url{http://www.ga.gov.au/scientific-topics/water/groundwater/groundwater-in-australia/fractured-rocks}, 2014.
\newblock [Online; accessed 21-February-2020].

\bibitem{GU2020109378}
{\sc Gu, Y., and Shen, J.}
\newblock Bound preserving and energy dissipative schemes for porous medium equation.
\newblock {\em Journal of Computational Physics 410\/} (2020), 109378.

\bibitem{Gustafson1994}
{\sc Gustafson, G., and Kr{\'a}sn{\'y}, J.}
\newblock Crystalline rock aquifers: Their occurrence, use and importance.
\newblock {\em Applied Hydrogeology 2}, 2 (1994), 64--75.

\bibitem{Harr2012groundwater}
{\sc Harr, M.}
\newblock {\em {Groundwater and Seepage}}.
\newblock {Dover Civil and Mechanical Engineering}. {Dover Publications}, 2012.

\bibitem{Izbash1931}
{\sc Izbash, S.~V.}
\newblock {\em O filtracii v krupnozernistom materiale.}
\newblock {I}zv. {N}auchno-{I}ssled. {I}nst. {G}idro-{T}ekh. ({N}.{I}.{L}{G}.), {L}eningrad 1, 1931 ({I}n~{R}ussian).

\bibitem{Jiang2021}
{\sc Jiang, Y.}
\newblock High order finite difference multi-resolution {WENO} method for nonlinear degenerate parabolic equations.
\newblock {\em Journal of Scientific Computing 86}, 1 (2021).

\bibitem{King1898}
{\sc King, F.}
\newblock Principles and conditions of the movements of ground water.
\newblock {\em Nineteenth Ann. Kept. U. S. Geol. Survey pt. 2}, 9-12 (1898), 209--215.

\bibitem{KIRALY1971255}
{\sc Kiràly, L.}
\newblock Groundwater flow in heterogeneous, anisotropic fractured media: A simple two-dimensional electric analog.
\newblock {\em Journal of Hydrology 12}, 3 (1971), 255--261.

\bibitem{Kroeber1884}
{\sc Kr{\"o}ber, C.}
\newblock Versuche \"uber die bewegung des wassers durch sandschichten.
\newblock {\em Zeitschr. des Vereines deutscher Ing. 28}, 31 and 32 (1884), 593--595 and 617 --619.

\bibitem{Lachassagne2021}
{\sc Lachassagne, P., Dewandel, B., and Wyns, R.}
\newblock Review: Hydrogeology of weathered crystalline/hard-rock aquifers---guidelines for the operational survey and management of their groundwater resources.
\newblock {\em Hydrogeology Journal 29}, 8 (Dec 2021), 2561--2594.

\bibitem{LapcevicNovakowskiSudicky}
{\sc Lapcevic, P.~A., Novakowski, K.~S., and Sudicky, E.~A.}
\newblock The interpretation of a tracer experiment conducted in a single fracture under conditions of natural groundwater flow.
\newblock {\em Water Resources Research 35}, 8 (1999), 2301--2312.

\bibitem{larsson1984ground}
{\sc Larsson, I., Unesco, and 8.6, I. H. P.~P.}
\newblock {\em Ground Water in Hard Rocks: Project 8.6 of the International Hydrological Programme}.
\newblock Studies and reports in hydrology. Unesco, 1984.

\bibitem{Leibenson1945}
{\sc Leibenson, L.~S.}
\newblock Turbulent movement of gas in a porous medium.
\newblock {\em Bull. Acad. Sci. USSR. S\'er. G\'eograph. G\'eophys. [Izvestia Akad. Nauk SSSR] 9\/} (1945), 3--6.
\newblock In Russian. Reprinted in Ref. \cite{Leibenson1953}, 499--502.

\bibitem{Leibenson1953}
{\sc Leibenson, L.~S.}
\newblock {\em Sobranie trudov, Chast' II: Podzemnaya gidrodinamika [Collected Works, Vol. II: Underground Hydrodynamics]}.
\newblock Izdat'elstvo Akademii Nauk S.S.S.R., Moscow, U.S.S.R., 1953.
\newblock In Russian.

\bibitem{Liu2011939}
{\sc Liu, Y., Shu, C.-W., and Zhang, M.}
\newblock High order finite difference {WENO} schemes for nonlinear degenerate parabolic equations.
\newblock {\em SIAM Journal on Scientific Computing 33}, 2 (2011), 939 – 965.

\bibitem{LUCAS2007295}
{\sc Lucas, Y., Panfilov, M., and Buès, M.}
\newblock High velocity flow through fractured and porous media: the role of flow non-periodicity.
\newblock {\em European Journal of Mechanics - B/Fluids 26}, 2 (2007), 295--303.

\bibitem{MacDonald2008}
{\sc Macdonald, A., Davies, J., and Calow, R.}
\newblock African hydrogeology and rural water supply.
\newblock {\em Applied Groundwater Studies in Africa\/} (2008), 127--148.

\bibitem{Marino1974}
{\sc Marino, M.}
\newblock Rise and decline of the water table induced by vertical recharge.
\newblock {\em Journal of Hydrology 23}, 3-4 (1974), 289--298.

\bibitem{Marechal2004}
{\sc Maréchal, J.~C., Dewandel, B., and Subrahmanyam, K.}
\newblock Use of hydraulic tests at different scales to characterize fracture network properties in the weathered-fractured layer of a hard rock aquifer.
\newblock {\em Water Resources Research 40}, 11 (2004).

\bibitem{Missbach1936}
{\sc Missbach, A.~A.}
\newblock Filtrovatelnost \v{c}e\v{r}en\'{y}ch a saturovan\'{y}ch \v{s}t'\'{a}v. {IV}. {P}\v{r}ezkou\v{s}en\'{i} vzorce van {G}ilse, ...
\newblock {\em Listy cukrov. 54}, 39 (1936), 361 -- 368.
\newblock In Czech.

\bibitem{OF}
{\sc {OpenFoam }}.
\newblock \url{ https://www.openfoam.com}.
\newblock [Online; accessed 02/28/2024 23:36].

\bibitem{Parlange2000339}
{\sc Parlange, J.-Y., Hogarth, W., Govindaraju, R., Parlange, M., and Lockington, D.}
\newblock On an exact analytical solution of the {B}oussinesq equation.
\newblock {\em Transport in Porous Media 39}, 3 (2000), 339 – 345.

\bibitem{Perrin2011}
{\sc Perrin, J., Ahmed, S., and Hunkeler, D.}
\newblock The effects of geological heterogeneities and piezometric fluctuations on groundwater flow and chemistry in a hard-rock aquifer, southern {I}ndia.
\newblock {\em Hydrogeology Journal 19}, 6 (2011), 1189--1201.

\bibitem{QUINN2020124384}
{\sc Quinn, P., Cherry, J., and Parker, B.}
\newblock Relationship between the critical reynolds number and aperture for flow through single fractures: Evidence from published laboratory studies.
\newblock {\em Journal of Hydrology 581\/} (2020), 124384.

\bibitem{SinghRai2D}
{\sc Rai, S., and Singh, R.}
\newblock Two-dimensional modelling of water table fluctuation in response to localised transient recharge.
\newblock {\em Journal of Hydrology 167}, 1 (1995), 167--174.

\bibitem{sarkar2004fluid}
{\sc Sarkar, S., Toksoz, M.~N., and Burns, D.~R.}
\newblock Fluid flow modeling in fractures.
\newblock Tech. rep., Massachusetts Institute of Technology. Earth Resources Laboratory, 2004.

\bibitem{Scheidegger1960}
{\sc Scheidegger, A.~E.}
\newblock {\em The Physics of Flow through Porous Media}.
\newblock The Macmillan company, New York, 1960.

\bibitem{ZekaiSen1995}
{\sc {\c S}en, Z.}
\newblock {\em Applied Hydrogeology for Scientists and Engineers}.
\newblock CRC Press, Boca Raton, 1995.

\bibitem{Shapiro2002}
{\sc Shapiro, A.~M.}
\newblock Fractured-rock aquifers understanding an increasingly important source of water.
\newblock \url{https://pubs.usgs.gov/publication/fs11202}, 2002.
\newblock [Online; accessed 02/28/2024 23:43].

\bibitem{SinghRai1980}
{\sc Singh, R., and Rai, S.}
\newblock On subsurface drainage of transient recharge.
\newblock {\em Journal of Hydrology 48}, 3-4 (1980), 303--311.

\bibitem{SinghRai1989}
{\sc Singh, R., and Rai, S.}
\newblock A solution of the nonlinear {B}oussinesq equation for phreatic flow using an integral balance approach.
\newblock {\em Journal of Hydrology 109}, 3-4 (1989), 319--323.

\bibitem{Singhal2008}
{\sc Singhal, B. B.~S.}
\newblock {\em Nature of Hard Rock Aquifers: Hydrogeological Uncertainties and Ambiguities}.
\newblock Springer Netherlands, Dordrecht, 2008, pp.~20--39.

\bibitem{Smreker1878}
{\sc Smreker, O.}
\newblock Entwicklung eines {G}esetzes f\"{u}r den {W}iderstand bei der{B}ewegung des {G}rundwassers.
\newblock {\em Zeitschr. des Vereines deutscher Ing. 22}, 4 and 5 (1878), 117--128 and 193--204.

\bibitem{Smreker1879}
{\sc Smreker, O.}
\newblock Das grundwasser und seine verwendung zu wasserversorgungen.
\newblock {\em Zeitschr. des Vereines deutscher Ing. 23}, 4 (1879), 347--362.

\bibitem{SoniIslamBasak1978}
{\sc Soni, J., Islam, N., and Basak, P.}
\newblock An experimental evaluation of non-darcian flow in porous media.
\newblock {\em Journal of Hydrology 38}, 3--4 (1978), 231--241.

\bibitem{StarkVolker}
{\sc Stark, K.~P., and Volker, R.~E.}
\newblock A study of some theoretical aspects of non-linear flow through porous materials.
\newblock Tech. Rep. in \emph{Res. Bulletin} No. 1, Dept. of Civil Engineering, University College of Townsville, Townsville, Australia, April 1967.

\bibitem{diplomka}
{\sc \v{S}vandov\'{a}, A.}
\newblock Simulation of fluid flow in porous media with an emphasis on fractured media and coarse-grained materials.
\newblock Master's thesis, University of West Bohemia, 2022.
\newblock Mentor:~P.~Girg.

\bibitem{Wright1992}
{\sc Wright, E.}
\newblock The hydrogeology of crystalline basement aquifers in {A}frica.
\newblock {\em Geological Society Special Publication 66\/} (1992), 1--27.

\bibitem{yan2018non}
{\sc Yan, X., Qian, J., Ma, L., Wang, M., and Hu, A.}
\newblock Non-fickian solute transport in a single fracture of marble parallel plate.
\newblock {\em Geofluids 2018\/} (2018).

\bibitem{Zhukovskii1889}
{\sc Zhukovskii, N.~E.}
\newblock Teoreticheskoe issledovanie o dvizhenii podpochvennykh vod.
\newblock {\em Zhurnal Russkogo fiziko\--khimicheskogo obshchestva 21}, 1 (1889).
\newblock [In Russian], reprinted in Ref.~\cite{Zhukovskii1937}, p. 9 -- 33.

\bibitem{Zhukovskii1937}
{\sc Zhukovskii, N.~E.}
\newblock {\em Polnoe sobranie sochinenii}, vol.~7.
\newblock Moscow\--Leningrad, 1937.
\newblock [Collected Papers], Russian with English Summary.

\bibitem{Zimmerman2004163}
{\sc Zimmerman, R.~W., Al-Yaarubi, A., Pain, C.~C., and Grattoni, C.~A.}
\newblock Non-linear regimes of fluid flow in rock fractures.
\newblock {\em International Journal of Rock Mechanics and Mining Sciences 41}, SUPPL. 1 (2004), 163 – 169.

\bibitem{Zunker1920}
{\sc Zunker, F.}
\newblock Das allgemeine {G}rund\-wasser\-fliess\-gesetz.
\newblock {\em Journal f\"{u}r {G}asbeleuchtung und {W}asserversorgung 63}, 21 (1920), 331--334, and 350.
\end{thebibliography}
\end{document}